\documentclass[11pt,twoside,a4paper]{article}

\usepackage{mathtools,bm,amssymb,amsmath,amssymb,mathrsfs,amsfonts}
\usepackage{graphicx}
\usepackage{float}
\usepackage{tikz}
\usepackage{tikz-cd}
\usetikzlibrary{shapes.geometric}

\usepackage[utf8]{inputenc}
\usepackage[T1]{fontenc}
\usepackage{amsmath,amssymb,dsfont}
\numberwithin{equation}{section}
\usepackage{microtype}
\usepackage{lmodern}
\usepackage{graphicx,tikz,pgfplots}
\usepackage{tikz-3dplot}
\usetikzlibrary{arrows,shapes,fadings,intersections,decorations.pathreplacing,positioning}
\graphicspath{{images/}}
\pgfplotsset{compat=newest}
\usepackage[hyperref,amsmath,thmmarks]{ntheorem}
\usepackage{aliascnt}
\usepackage[a4paper,centering,bindingoffset=0cm,marginpar=2cm,margin=2.5cm]{geometry}
\usepackage[pagestyles]{titlesec}
\usepackage[font=footnotesize,format=plain,labelfont=sc,textfont=sl,width=0.75\textwidth,labelsep=period]{caption}
\usepackage{comment}
\usepackage{siunitx}
\usepackage{mathtools}

\usetikzlibrary{quotes}
\usetikzlibrary{angles}
\usetikzlibrary{babel}

\usetikzlibrary{decorations.pathreplacing}
\usetikzlibrary{decorations.markings}

\usepackage{mathabx}
\tikzset{
    >=stealth',
    punkt/.style={
           rectangle,
           rounded corners,
           draw=black, very thick,
           text width=6.5em,
           minimum height=2em,
           text centered},
    pil/.style={
           ->,
           thick,
           shorten <=2pt,
           shorten >=2pt,}
}
\tikzstyle{block} = [rectangle, rounded corners, minimum width= 3cm, minimum height=1cm, text centered, draw=black, fill=blue!20,]

\tikzstyle{invisbleblock} = [minimum width= 3em, minimum height=1cm, text centered, ]

\tikzstyle{decision} = [diamond, minimum width=3cm, minimum height=1cm, text centered, draw=black, fill=orange!20]

\tikzstyle{arrow} = [thick,->,>=stealth]
\tikzstyle{line} = [thick,-]

\usepackage{dsfont,bbm}

\usepackage[boxed]{algorithm2e}

\usepackage{enumerate}

\usepackage{subcaption}
\usepackage{multirow}

\usepackage[backend=biber,maxnames=10,backref=true,hyperref=true,giveninits=true,safeinputenc]{biblatex}
\DefineBibliographyStrings{english}{%
	backrefpage = {cited on page},
	backrefpages = {cited on pages},
}

\title{A new inversion scheme for elastic diffraction tomography}
\author{Bochra Mejri$^1$\\
	{\footnotesize\href{mailto:bochra.mejri@ricam.oeaw.ac.at}{bochra.mejri@ricam.oeaw.ac.at}}
	\and Otmar Scherzer$^{1,2,3}$\\
	{\footnotesize\href{mailto:otmar.scherzer@univie.ac.at}{otmar.scherzer@univie.ac.at}}
}

\date{\today}

\DeclareFieldFormat[report]{title}{``#1''}
\DeclareFieldFormat[book]{title}{``#1''}

\titleformat{\section}[block]{\large\sc\filcenter}{\thesection.}{0.5ex}{}[]
\titleformat{\subsection}[runin]{\bf}{\thesubsection.}{0.5ex}{}[.]
\titleformat{\subsubsection}[runin]{\bf}{\thesubsubsection.}{0.5ex}{}[.] 

\usepackage[pdftex,colorlinks=true,linkcolor=blue,citecolor=green!50!black,urlcolor=blue,bookmarks=true,bookmarksnumbered=true,pdfusetitle]{hyperref}
\hypersetup
{
    pdfauthor={B.~Mejri and O.~Scherzer},
}

\newpagestyle{headers}
{
	\headrule
	\sethead[\footnotesize\thepage][\footnotesize\sc B.~Mejri and O.~Scherzer][]{}{\footnotesize\sc A new inversion scheme for elastic diffraction tomography
  }{\footnotesize\thepage}
	\setfoot{}{}{}
}
\newtheorem{lemma}{Lemma}[section]

\newaliascnt{proposition}{lemma}

\aliascntresetthe{proposition}

\newaliascnt{corollary}{lemma}
\newtheorem{corollary}[corollary]{Corollary}
\aliascntresetthe{corollary}

\newaliascnt{theorem}{lemma}
\newtheorem{theorem}[theorem]{Theorem}
\aliascntresetthe{theorem}

\theorembodyfont{\normalfont}
\newaliascnt{definition}{lemma}

\aliascntresetthe{definition}

\newaliascnt{assumption}{lemma}

\aliascntresetthe{assumption}

\newaliascnt{notation}{lemma}

\aliascntresetthe{notation}

\newaliascnt{example}{lemma}
\newtheorem{example}[example]{Example}
\aliascntresetthe{example}

\newaliascnt{experiment}{lemma}

\aliascntresetthe{experiment}

\newaliascnt{remark}{lemma}
\newtheorem{remark}[remark]{Remark}
\aliascntresetthe{remark}

\theoremstyle{nonumberplain}
\theoremseparator{:}
\theoremheaderfont{\normalfont\itshape}

\theoremsymbol{\ensuremath{\square}}
\newtheorem{proof}{Proof}

\newcommand{\vecu}{\bm{u}}
\newcommand{\vecut}{\bm{u}^\text{tot}}
\newcommand{\vecui}{\bm{u}^\text{inc}}
\newcommand{\vecus}{\bm{u}^\text{sca}}
\newcommand{\vecub}{\bm{u}^\text{Born}}
\newcommand{\vecutt}{\bm{u}_t}

\newcommand{\vecn}{\bm{n}}
\newcommand{\vecr}{\bm{r}}
\newcommand{\vece}{\bm{e}}
\newcommand{\vecq}{\bm{q}}
\newcommand{\vecf}{\bm{f}}
\newcommand{\vecg}{\bm{g}}
\newcommand{\vech}{\bm{h}}

\newcommand{\sigmabo}{\boldsymbol{\sigma}}
\newcommand{\sigmabot}{\boldsymbol{\sigma}^\text{tot}}
\newcommand{\sigmaboi}{\boldsymbol{\sigma}^\text{inc}}
\newcommand{\sigmabos}{\boldsymbol{\sigma}^\text{sca}}
\newcommand{\sigmabob}{\boldsymbol{\sigma}^\text{Born}}
\newcommand{\sigmabott}{\boldsymbol{\sigma}_t}

\newcommand{\varepsilonbo}{\boldsymbol{\varepsilon}}

\newcommand{\nablabo}{\boldsymbol{\nabla}}

\newcommand{\varphibo}{\boldsymbol{\varphi}}
\newcommand{\xibo}{\boldsymbol{\xi}}
\newcommand{\vecC}{\bm{\mathcal{C}}}
\newcommand{\vecG}{\bm{\mathcal{G}}}

\newcommand{\x}{\bm{x}}
\newcommand{\y}{\bm{y}}

\newcommand{\dy}{\,\mbox{d}\bm{y}}

\def\N{\mathbb{N}}
\def\R{\mathbb{R}}
\def\C{\mathbb{C}}
\def\S{\mathbb{S}}

\begin{document}

\maketitle
\thispagestyle{empty}
\begin{center}
	\parbox[t]{15em}{\footnotesize
		\hspace*{-1ex}$^1$Johann Radon Institute 
		for Computational and Applied Mathematics\\ (RICAM)\\
		Altenbergerstraße 69\\
		A-4040 Linz, Austria} 
	\hspace{1em}
	\parbox[t]{11em}{\footnotesize
		\hspace*{-1ex}$^2$Faculty of Mathematics\\
		University of Vienna\\
		Oskar-Morgenstern-Platz 1\\
		A-1090 Vienna, Austria}
	\hspace{1em}
	\parbox[t]{11em}{\footnotesize
		\hspace*{-1ex}$^3$Christian Doppler Laboratory
		 for Mathematical Modeling and Simulation
		  of Next Generations of Ultrasound Devices (MaMSi)\\
		Institute of Mathematics\\
		Oskar-Morgenstern-Platz 1\\
		A-1090 Vienna, Austria}
\end{center}

\begin{abstract}
We consider the problem of elastic diffraction tomography, which consists in reconstructing elastic properties (i.e. mass density and elastic Lamé parameters) of a weakly scattering medium from full-field data of scattered waves outside the medium. Elastic diffraction tomography refers to the elastic inverse scattering problem after linearization using a first-order Born approximation. In this paper, we prove the Fourier diffraction theorem, which relates the 2D Fourier transform of scattered waves with the Fourier transform of the scatterer in the 3D spatial Fourier domain. Elastic wave mode separation is performed, which decomposes a wave into four modes. A new two-step inversion process is developed, providing information on the modes first and secondly on the elastic parameters. Finally, we discuss reconstructions with plane wave excitation experiments for different tomographic setups and with different plane wave excitation frequencies, respectively.
\end{abstract}



\section{Introduction}
\label{introduction}
Diffraction tomography (acoustic, elastic) continues to be a thriving mathematical research interest. It has widespread applications in nondestructive testing, geophysical explorations and medical diagnostics~\cite{LanMayMar06,LanMayMar10,Dev84}. Acoustic diffraction tomography is an inverse scattering technique aiming to imaging refractive index distribution of a \emph{weakly scattering medium} with probing waves~\cite{KakSla01,NatWub01}. More generally, \emph{Elastic diffraction tomography} seeks to retrieve quantitative information on elastic parameters (i.e. mass density and elastic Lamé parameters) from measurement of scattered waves recorded outside the medium. 

In an elastic medium, there are two modes of wave propagation, \emph{pressure waves} and \emph{shear waves}, while in an acoustic medium, only pressure waves propagate. Mechanical waves are characterised by different wave speed (i.e. pressure wave propagate with a velocity faster than shear wave for the same frequency) and polarization directions (i.e. pressure polarization is parallel to the wave direction, while shear polarization is in the perpendicular plane). Pressure waves may convert to shear waves through refraction or reflection at a boundary and vice versa. Then, the scattered waves contain a superposition of both wave modes at the boundary of the medium discontinuities. 

The difficulties affecting the study of elastic scattering problems are twofold: wave mode separation and multi-parameter inversion process (i.e. three elastic parameters are to be quantified). Therefore, one can find in the literature various approaches for elastic parameter reconstruction based on certain assumptions. The inverse problem can be linearized using Born~\cite{GubDomKruHub77,BlaBurHopWom87,BeyBur90,LanMayMar06} or Rytov approximation~\cite{GelVirGra07}. Afterwards, the developed algorithms are based on system reduction (e.g. scalar elastodynamic equation~\cite{BlaBurHopWom87}, one perturbed parameter~\cite{LanMayMar06}) or/and approximations of the measurements (e.g. near or far-field approximation~\cite{BlaBurHopWom87} or considering singular part~\cite{BeyBur90}). The approaches therein are based on two step reconstruction process. A first estimation on the three-parameters (i.e. mass density and elastic Lamé parameters) is provided from measurements on three different positions of source/receiver. Secondly, a higher resolution image is obtained by making a full rotation of the object~\cite{BlaBurHopWom87} or by sweeping the frequency bandwidth~\cite{BeyBur90}.

For \emph{acoustic diffraction}, Kirisits {\it et al.}~\cite{KirQueRitSchSet21} (i.e. scalar theory) established a detailed mathematical analysis of the Fourier diffraction theorem and a backprojection formula. Besides, they discussed reconstruction formulas for different tomographic setups~\cite{KirQueRitSchSet21} and with plane wave excitations of different frequencies~\cite{FauKirQueSchSet21_report}. The main objective of this paper is twofold: 
\begin{enumerate}
\item {\bf From acoustic to elastic:} First, we extend the studies from the acoustic case~\cite{KirQueRitSchSet21} to the elastic case. This problem requires more subtle analysis due to the coupling between the pressure and shear waves. Furthermore, the inversion process is not straightforward due to the singularities of the partial Fourier transform of Green's tensor.
\item {\bf General framework}: Secondly, we consider a general framework of elastodynamic system integrating mode superposition. In other words, we consider a full 3D vector elasticity system characterised by the mass density and the Lamé parameters. To our best knowledge, the present paper is a first attempt for the parallel reconstruction of the three-elastic parameters from a full-field data of scattered waves (i.e. no approximation on the measurements). 
\end{enumerate}

{\bf Contributions and outline.}
Our approach relies on the Fourier diffraction theorem related to the 2D Fourier transform of full-field measurements of diffracted waves and to the Fourier transform of the scatterer (i.e. referring to as scattering potential), which is mapped into the 3D spatial Fourier space (i.e. k-space). The decomposition provides in total four elastic wave modes PP, PS, SP, SS (i.e. P: pressure, S: shear) containing complementary information on the data, form which the elastic parameters could be reconstructed. Here, we develop a new multi-parameter backprojection inversion scheme for the multi-modes. Finally, we derive k-space coverage for each mode by varying the excitation direction or by sweeping the frequency bandwidth.  

This article is outlined as follows. The mathematical basics of elastic scattering are reviewed in subsection~\ref{inhomogeneous_elastic_wave}. The elastic diffraction model, obtained by linearization using a first-order Born approximation, is given in subsection~\ref{born_approximation} and the Helmholtz's decomposition for the incident waves is considered in subsection~\ref{helmholtzs_decomposition}. The Fourier diffraction theorem for transmission and reflection acquisitions experiments is presented in subsection~\ref{fourier_diffraction_theorem} and its proof is given in Appendix~\ref{appendix}. In subsection~\ref{wave_mode_separation}, wave mode separation is performed using the properties of shear and pressure waves and specific filters. Section~\ref{wave_field_inversion} presents a new multi-parameter inversion algorithm, which operates in the Fourier domain. Different coverages of k-space are obtained in terms of angular diversity in subsection~\ref{angular_diversity} and in terms of frequency diversity in subsection~\ref{frequency_diversity}.

\section{Elastic waves}
\label{elastic_waves}
\subsection{Homogeneous elastic wave}
\label{homogeneous_elastic_wave}
The considered domain is an infinite, isotropic elastic body characterized by the mass density $\rho\in L^{\infty}(\R^3,\R_+)$ and the fourth-order elasticity tensor $\vecC$ with components $\mathcal{C}_{ijkl}\in L^{\infty}(\R^3,\R)$ defined at a point $\x\in\R^3$ as
\[
	\mathcal{C}_{ijkl}(\x)=\lambda(\x)\delta_{ij}\delta_{kl}+\mu(\x)(\delta_{ki}\delta_{lj}+\delta_{li}\delta_{kj}),\quad\quad\text{for }i,j,k,l=1,2,3
\]
where the Lamé parameters $\mu,\lambda\in  L^{\infty}(\R^3,\R)$ with $\mu>0$ and $\lambda+2\mu>0$, while $\delta_{ij}$ is the Kronecker’s delta. 
    
The mechanical problem considered is represented by the following system that expresses the propagation of time-harmonic elastic waves in the reference background. 
The incident displacement field $\vecui(\x,\omega)$ satisfies  
    \begin{equation}\label{inc_sys}
    \left\{\begin{aligned}
    & \sigmaboi(\x,\omega)=\vecC^0(\x):\varepsilonbo\left[\vecui\right](\x,\omega)                       \\
    & \nablabo\cdot\sigmaboi(\x,\omega)+\omega^2\rho^0(\x)\vecui(\x,\omega)= \bm{0}              
    \end{aligned} \right.
    \end{equation}
with linearized strain tensor $\varepsilonbo\left[\vecu\right](\x,\omega)=\frac12(\nablabo\vecu(\x,\omega)+\nablabo\vecu^\top(\x,\omega))$ and stress field $\sigmabo(\x,\omega)$ at an angular frequency $\omega>0$, with the symbols $"\cdot"$ and $":"$ respectively denote simple and double tensor inner products.

\subsection{Inhomogeneous elastic wave}
\label{inhomogeneous_elastic_wave}
Considering the scattering by a weak heterogeneity of an elastic plane wave propagating in an infinite medium (i.e. reference background) defined by a density $\rho^0$ and an elastic stiffness $\vecC^0(\mu^0,\lambda^0)$. The heterogeneity is characterized by perturbed parameters: a density $\delta\rho$ and an elastic stiffness $\delta\vecC(\delta\mu,\delta\lambda)$. One assumes that the object is included in an open ball $\mathcal{B}_{r_s}\subset\R^3$ with midpoint $\bm{0}$ and radius $r_s$. 

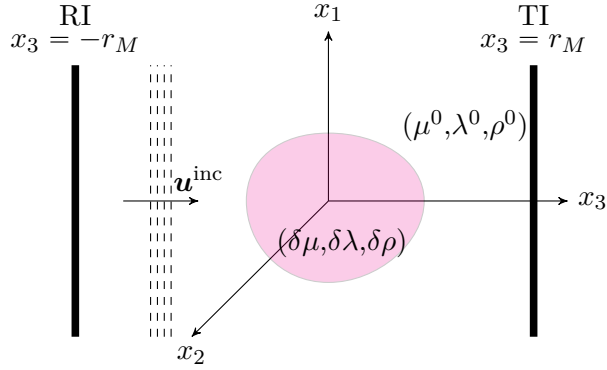
\begin{figure}[h!]  
\centering 
\begin{tikzpicture}[scale=.9]
	\draw[->] (0,0) -- (0,2.5) node[anchor = south]{$x_1$};
	\draw[->] (0,0) -- (-2,-2) node[anchor = north]{$x_2$};
	\draw[->] (0,0) -- (3.5,0) node[anchor = west]{$x_3$};
        \draw [black,fill = magenta,opacity=.2] plot[smooth cycle,tension=1] coordinates {(0,-1.2)(-1.2,0)(0,1.)(1.4,0)};
        \draw [thick] (0.2,-.25) node[anchor = north]{$(\delta\mu$,$\delta\lambda$,$\delta\rho)$};
        \draw [thick] (2,1.5) node[anchor = north]{$(\mu^0$,$\lambda^0$,$\rho^0)$};
        \draw [line width=1mm] (3,-2) -- (3,2) node[anchor = south] {$x_3=r_M$} node[above =10pt] {TI};
        \draw [line width=1mm] (-3.7,-2) -- (-3.7,2) node[anchor = south] {$x_3=-r_M$}node[above =10pt] {RI};
        \draw[->] (-3,0) -- (-1.9,0) node[anchor = south]{$\vecui$};
        \draw [dashed] (-2.6,-2) -- (-2.6,2);
        \draw [dashed] (-2.5,-2) -- (-2.5,2);
        \draw [dashed] (-2.4,-2) -- (-2.4,2);
        \draw [dashed] (-2.3,-2) -- (-2.3,2);
\end{tikzpicture}
\caption{Elastic scattering situation for a transmission imaging (TI) and a reflection imaging (RI).}
\label{four_diff_figa}
\end{figure}

The elastic scattering problem is defined as follows: Firstly, one generates incident waves propagating in the direction of the heterogeneity which selectively scatters part of the incident power into various directions. Then, one measures the incident waves scattered into a given direction at a distance from the object $x_3=r_M>r_s$ for transmission imaging and $x_3=-r_M<-r_s$ for reflection imaging, see Figure~\ref{four_diff_figa}.    

Therefore, there is a total wave $\vecut$ (i.e. total displacement field) written as the sum of an incident wave (i.e. $\vecui$ incident displacement field) and a scattered wave (i.e. $\vecus$ scattered displacement field) given by 
\[
\vecut=\vecui+\vecus
\]
Besides, the total displacement field $\vecut(\x,\omega)$ satisfies  
   \begin{equation}\label{tot_sys}
    \left\{\begin{aligned}
    & \sigmabot(\x,\omega)=\vecC(\x):\varepsilonbo\left[\vecut\right](\x,\omega)                       \\
    & \nablabo\cdot\sigmabot(\x,\omega)+\omega^2\rho(\x)\vecut(\x,\omega)= \bm{0}              
    \end{aligned} \right.
    \end{equation}
where $\vecC(\x)=\vecC^0(\x)+\delta\vecC(\x)$ is the total elasticity tensor and $\rho(\x)=\rho^0(\x)+\delta\rho(\x)$ is the total mass density. 

\subsection{Born approximation}
\label{born_approximation}
To linearize the {\it{nonlinear}} inverse problem, one can use the Kirchhoff approximation for strong scatterers~\cite{LanMayMar06} or the Born approximation for weak scatterers~\cite{GubDomKruHub77}, the one being considered herein. Then, one gets by inserting Equation~\ref{inc_sys} in Equation~\ref{tot_sys} with the fact that $\vecut=\vecui+\vecus$, $\vecC(\x)=\vecC^0(\x)+\delta\vecC(\x)$ and $\rho(\x)=\rho^0(\x)+\delta\rho(\x)$
   \begin{equation}\label{sca_sys}
    \left\{\begin{aligned}
    & \sigmabos(\x,\omega)=\vecC^0(\x):\varepsilonbo\left[\vecus\right](\x,\omega)     &                 \\
    & \nablabo\cdot\sigmabos(\x,\omega)+\omega^2\rho^0(\x)\vecus(\x,\omega)=&-\nablabo\cdot\left(\delta\vecC(\x):\varepsilonbo\left[\vecui\right](\x,\omega)\right)-\omega^2\delta\rho(\x)\vecui(\x,\omega)\\
    &&-\nablabo\cdot\left(\delta\vecC(\x):\varepsilonbo\left[\vecus\right](\x,\omega)\right)-\omega^2\delta\rho(\x)\vecus(\x,\omega)           
    \end{aligned} \right.
    \end{equation}
Neglecting $\vecus$ in the right-hand side of Equation~\ref{sca_sys}  (assuming that $\vecus$ is small compared to $\vecui$), one obtains the first-order {\it Born approximation} as follows 
  \begin{equation}\label{born_app}
    \left\{\begin{aligned}
    & \sigmabob(\x,\omega)=\vecC^0(\x):\varepsilonbo\left[\vecub\right](\x,\omega)                      \\
    & \nablabo\cdot\sigmabob(\x,\omega)+\omega^2\rho^0(\x)\vecub(\x,\omega) = -\nablabo\cdot\left(\delta\vecC(\x):\varepsilonbo\left[\vecui\right](\x,\omega)\right)-\omega^2\delta\rho(\x)\vecui(\x,\omega)           
    \end{aligned} \right.
    \end{equation}
The vector function 
\[
\vecf(\x,\omega)=-\nablabo\cdot\left(\delta\vecC(\x):\varepsilonbo\left[\vecui\right](\x,\omega)\right)-\omega^2\delta\rho(\x)\vecui(\x,\omega)           
\]  
is referred as the {\it{scattering potential}} and will be the quantity which we set out to reconstruct from the measurements of the scattered waves. By construction, we have 
\[
\text{supp}(\vecf)\subseteq\mathcal{B}_{r_s}\subset(-r_s,r_s)^3
\]    		
For simplicity of notation, we set $\vecub=\vecu$ from now on.  
\begin{remark}\label{born_approx_remark}
The Born approximation is valid under the condition that the total phase fluctuation caused by scattering must be less than one radian, for more details see~\cite{Wu89}. 
\end{remark} 

\subsection{Helmholtz's decomposition}
\label{helmholtzs_decomposition}
The displacement $\vecu$, satisfying {\it{homogeneous elastodynamic}} system~\ref{born_app} with $\vecf=\bm{0}$, can be written as the sum of a shear part $\vecu_s$ (i.e. transverse component) and a pressure part $\vecu_p$ (i.e. longitudinal component) using the {\it{Helmholtz decomposition}} as follows 
\[
\vecu=\vecu_s+\vecu_p
\] 
where $\vecu_s$ is the divergence-free part (i.e. $\nabla\cdot\vecu_s=0$) and $\vecu_p$ is the curl-free part (i.e. $\nabla\times\vecu_p=\bm{0}$) of the displacement $\vecu$. Denoting by $k_s$ and $k_p$ the wavenumber of shear waves and pressure waves, respectively, defined as follows  
\[
k_s^2=\frac{\omega^2\rho^0}{\mu^0}\quad\text{and}\quad k_p^2=\frac{\omega^2\rho^0}{\lambda^0+2\mu^0}
\]
Besides, $\vecu_s$ and $\vecu_p$ are solutions to the {\it{vectoriel homogeneous Helmholtz}} equations $\Delta\vecu_s(\x)+k_s^2\vecu_s(\x)=\bm{0}$ and $\Delta\vecu_p(\x)+k_p^2\vecu_p(\x)=\bm{0}$, respectively (for more details see Section 5.16 in~\cite{EriSuhCha78}).

If $\vecu$ is a solution defined in an infinite elastic medium satisfying the following outgoing Kupradze radiation conditions 
\begin{equation}\label{eq:kuprad}
\lim_{r\to\infty}\max_{\|\vecr\|=r}\|\vecr\|\left(\partial_{\vecr}-ik_s\right)\vecu_s(\vecr)=0\quad\text{and}\quad\lim_{r\to\infty}\max_{\|\vecr\|=r}\|\vecr\|\left(\partial_{\vecr}-ik_p\right)\vecu_p(\vecr)=0
\end{equation}
where $\partial_{\vecr}$ denotes the directional derivative. Then, $\vecu=\bm{0}$. Physically speaking, $\vecu$ is an outgoing wave.
\begin{remark}\label{homo_elas_wave_remark} 
The pressure velocity is generally bigger than the shear velocity (i.e. $\lambda>0$) so the subscript $p$ stands for primary and $s$ for secondary waves indicating the earlier arrival time of P-waves against S-waves. P-waves are curl-free pressure waves and S-waves are purely solenoidal shear waves.
\end{remark}
According to Helmholtz's decomposition, one can consider the incident displacement field $\vecui$ as a trial solution for plane waves propagating in the $\vece_3=(0,0,1)^\top$ wave propagation direction written as follows   
\begin{equation}\label{inc_equ}
 \vecui(\x,\omega)=  \vecui_s(\x,\omega)+ \vecui_p(\x,\omega)
 \end{equation}
where the shear incident plane wave $\vecui_s(\x,\omega)=\mathtt{A}_s e^{ik_sx_3}$ with the shear amplitude vector $\mathtt{A}_s$ orthogonal to the wave propagation direction, while the pressure incident plane wave $\vecui_p(\x,\omega)=\mathtt{A}_p e^{ik_px_3}$ with the pressure amplitude vector $\mathtt{A}_p$ parallel to the vector $\vece_3$. In other words, the components $a_{s,3}=a_{p,1}=a_{p,2}$ are generally equal to zero and $a_{s,1},a_{s,2},a_{p,3}$ are not equal to zero.

The scattered wave $\vecu$ satisfies Equation~\ref{born_app} under the validity of the Born approximation (see Remark~\ref{born_approx_remark}). By exploiting the specific form of the incident wave $\vecui$, Equation~\ref{inc_equ}, the scattered displacement field $\vecu$ will be the solution of 
\begin{equation}\label{eq:bornapp}
    \left\{\begin{aligned}
    & \sigmabo(\x,\omega)=\vecC^0(\x):\varepsilonbo\left[\vecu\right](\x,\omega)                      \\
    & \nablabo\cdot\sigmabo(\x,\omega)+\omega^2\rho^0(\x)\vecu(\x,\omega) = \vecf_s(\x,\omega)e^{ik_sx_3} + \vecf_p(\x,\omega)e^{ik_px_3}              
    \end{aligned} \right.
\end{equation}
with the shear scattering potential $\vecf_s$ is given by 
\begin{equation}\label{fs_eq}
\vecf_s(\x,\omega):=k_s^2\delta\mu(\x)\mathtt{A}_s-ik_s\nablabo\cdot\sigmaboi_s(\x)-\omega^2\delta\rho(\x)\mathtt{A}_s
\end{equation}
where
\[
   \sigmaboi_s(\x)=\left(
    \begin{array}{ccc} 
    0            & 0    & \delta\mu(\x)a_{s,1}    \\
    0    & 0           & \delta\mu(\x)a_{s,2}     \\
    \delta\mu(\x)a_{s,1}     & \delta\mu(\x)a_{s,2}   & 0            \\
    \end{array}
    \right)
\]
and the pressure scattering potential $\vecf_p$ is given by 
\begin{equation}\label{fp_eq}
\vecf_p(\x,\omega):=k_p^2\left(\delta\lambda(\x)+2\delta\mu(\x)\right)\mathtt{A}_p-ik_p\nablabo\cdot\sigmaboi_p(\x)-\omega^2\delta\rho(\x)\mathtt{A}_p
\end{equation}
where
\[
    \sigmaboi_p(\x)=\left(
    \begin{array}{ccc} 
    \delta\lambda(\x)a_{p,3}            & 0    & 0    \\
    0    & \delta\lambda(\x)a_{p,3}           & 0     \\
    0     & 0   & (\delta\lambda(\x)+2\delta\mu(\x))a_{p,3}            \\
    \end{array}
    \right)
\]

\section{Fourier diffraction theorem}
To formulate the Fourier diffraction Theorem~\ref{four_diff_theo}, one needs the following notation
\begin{equation}
    \kappa_{\alpha}=\kappa_{\alpha}(\xibo):=\left\{
    \begin{aligned}
    &  \sqrt{k_{\alpha}^2-\xi_1^2-\xi_2^2},&       \xi_1^2+\xi_2^2\leqslant k_{\alpha}^2        &     \\
    & i\sqrt{\xi_1^2+\xi_2^2-k_{\alpha}^2},&      \xi_1^2+\xi_2^2>k_{\alpha}^2                     & \quad\mbox{for }{\alpha}=s,\,p
    \end{aligned} 
    \right. 
\end{equation}

Furthermore, $H_{x_3}:\R^3\to\R$ denotes the Heaviside function in the third coordinate centred at $x_3$, that is, 
\[
    H_{x_3}(s_1,s_2,s_3):=\left\{
    \begin{aligned}
    & 0 &       \mbox{if } s_3<x_3            \\
    & 1 &       \mbox{otherwise}.  
    \end{aligned} 
    \right.
\]
The partial Fourier transform $\mathcal{F}_{1,2}\vecu$ of a generic scalar, vector or tensor field $\vecu$ is defined by 
\[
\mathcal{F}_{1,2}\vecu(\xibo,x_3)=(2\pi)^{-1}\int_{\R^2}e^{-i(\xi_1x_1+\xi_2x_2)}\vecu(x_1,x_2,x_3)\,\mbox{d}x_1\mbox{d}x_2,
\]
where the spatial Fourier vector $\xibo:=(\xi_1,\xi_2)$ are the transformed coordinates associated with $(x_1,x_2)$. Finally, we denote by $\it{D}'(\R^3)$ the space of distribution and by $\it{S}'(\R^3)$ the space of tempered distributions.

\subsection{Fourier diffraction theorem}
\label{fourier_diffraction_theorem}
Here, one presents the Fourier diffraction theorem where the proof is postponed to Appendix~\ref{appendix}. 
\begin{theorem}
\label{four_diff_theo}
Let $\vecg\in\left[L^p(\R^3)\right]^3$, $p>1$, with $supp(\vecg)\subset\mathcal{B}_r$ for $r>0$. Suppose that the vector field $\vecu$ is the solution of 
\begin{equation}\label{eq32}
    \left\{\begin{aligned}
    & \sigmabo(\x,\omega)=\vecC^0(\x):\varepsilonbo\left[\vecu\right](\x,\omega)                      \\
    & \nablabo\cdot\sigmabo(\x,\omega)+\omega^2\rho^0(\x)\vecu(\x,\omega) = -\vecg(\x,\omega)           
    \end{aligned} \right.
\end{equation}
which satisfies the outgoing Kupradze radiation conditions, Equation~\ref{eq:kuprad}. Then, one can identify $\mathcal{F}_{1,2}\vecu\in\left[\mathcal{D}'(\R^3)\right]^3$ with the following locally integrable function
\begin{equation}\label{eq:fourdiff}
\begin{split}
 \mathcal{F}_{1,2}\vecu(\xibo,x_3,\omega)&=\sqrt{\frac{\pi}{2}}[\hat{\vecG}_s\left\{e^{i\kappa_sx_3}\mathcal{F}\left(\left(1-H_{x_3}\right)\vecg\right)(\xibo,\kappa_s,\omega)+e^{-i\kappa_sx_3}\mathcal{F}\left(H_{x_3}\vecg\right)(\xibo,-\kappa_s,\omega)\right\}
                                                         \\&+\hat{\vecG}_p\left\{e^{i\kappa_px_3}\mathcal{F}\left(\left(1-H_{x_3}\right)\vecg\right)(\xibo,\kappa_p,\omega)+e^{-i\kappa_px_3}\mathcal{F}\left(H_{x_3}\vecg\right)(\xibo,-\kappa_p,\omega)\right\}]
\end{split}
\end{equation}
for all $\xibo\in\R^2$ satisfying $\xi_1^2+\xi_2^2\neq k_s^2$ and $\xi_1^2+\xi_2^2\neq k_p^2$, where 
\begin{equation}\label{eq:ft_green_tensor}
\hat{\vecG}_s=\frac{i}{\mu k_s^2\kappa_s}\left[k_s^2\bm{I}_2+\vecq_s\otimes\vecq_s\right],\quad\quad\hat{\vecG}_p=\frac{-i}{\mu k_s^2\kappa_p}\vecq_p\otimes\vecq_p
\end{equation}
and having set the propagation vectors 
\begin{equation}\label{eq:prop_vect}
\vecq_s(\xibo):=-i\xibo'+i\kappa_s(\xibo)\text{sign}(x_3)\vece_3,\quad\quad\vecq_p(\xibo):=-i\xibo'+i\kappa_p(\xibo)\text{sign}(x_3)\vece_3
\end{equation}
with $\xibo'=(\xibo,0)$ and the symbol $\otimes$ denoting the dyadic product between tensors.
\end{theorem}
\begin{remark}
\begin{itemize}
\item The regularity choice for the second member $\vecg$ of system~\ref{eq32} is only for numerical purposes. 
\item In the case of evanescent waves (i.e. $\xi_1^2+\xi_2^2>k_{\alpha}^2$), one has to consider the analytic continuation of $\mathcal{F}\left((1-H_{x_3})\vecg\right)$ and $\mathcal{F}\left(H_{x_3}\vecg\right)$ to $\C^3$. 
\end{itemize}
\end{remark}
Fourier diffraction theorem~\ref{four_diff_theo} implies the following result for our setting. 
\begin{corollary}
\label{four_diff_coro}
Assume that the scattered displacement field $\vecu$ is the solution of system~\ref{eq:bornapp} satisfying the outgoing Kupradze radiation conditions, Equation~\ref{eq:kuprad}, with the shear scattering potential $\vecf_s\in\left[L^p(\R^3)\right]^3$, $p>1$, Equation~\ref{fs_eq}, with $supp(\vecf_s)\subset\mathcal{B}_{r_s}$ for $0<r_s<r_M$ and the pressure scattering potential $\vecf_p\in\left[L^p(\R^3)\right]^3$, $p>1$, Equation~\ref{fp_eq}, with $supp(\vecf_p)\subset\mathcal{B}_{r_s}$ for $0<r_s<r_M$. Then 
\begin{equation}\label{four_diff_eq}
\begin{split}
 \mathcal{F}_{1,2}\vecu(\xibo,\pm r_M,\omega)&=\sqrt{\frac{\pi}{2}}[\hat{\vecG}_se^{i\kappa_sr_M}\left\{\mathcal{F}\vecf_s(\xibo,\pm\kappa_s-k_s,\omega)+\mathcal{F}\vecf_p(\xibo,\pm\kappa_s-k_p,\omega)\right\}
                                                                \\&+\hat{\vecG}_pe^{i\kappa_pr_M}\left\{\mathcal{F}\vecf_s(\xibo,\pm\kappa_p-k_s,\omega)+\mathcal{F}\vecf_p(\xibo,\pm\kappa_p-k_p,\omega)\right\}]
 \end{split}
\end{equation}
for all $\xibo\in\R^2$ satisfying $\xi_1^2+\xi_2^2\neq k_s^2$ and $\xi_1^2+\xi_2^2\neq k_p^2$.
\end{corollary}
\begin{proof}
According to Equation~\ref{eq:fourdiff} with $\vecg=\vecf_se^{ik_sx_3}+\vecf_pe^{ik_px_3}$, one obtains $H_{r_M}\vecf_{\alpha}=\bm{0}$ and $(1-H_{r_M})\vecf_{\alpha}=\vecf_{\alpha}$, for transmission imaging (i.e. $x_3=r_M$). Similarly, for reflection imaging (i.e. $x_3=-r_M$), one gets $H_{r_M}\vecf_{\alpha}=\vecf_{\alpha}$ and $(1-H_{r_M})\vecf_{\alpha}=\bm{0}$.
\end{proof}

\subsection{Wave mode separation}
\label{wave_mode_separation}
Our aim, in this work, is to reconstruct the scattering potentials $\vecf_s$ and $\vecf_p$, defined in a volume $\mathcal{Y}$, from the scattered displacement data. To do so, we call the inverse spatial Fourier transform for a function $g$ as follows 
\[
g(\vecr)\approx(2\pi)^{-\frac32}\int_{\mathcal{Y}}e^{i\y\cdot\vecr}\mathcal{F}g(\y)\dy.
\]
It is a challenging task to invert the formula~\ref{four_diff_eq} due to the four terms defined in different k-spaces. At this point, one needs to consider term separations for the inversion process, in other words to perform a wave mode separation. More precisely, there is different modes of incident and scattered waves defined by shear or pressure parts. Then, the first step is to separate the shear part, see Equation~\ref{shear_four_diff_eq}, accordingly the pressure part, see Equation~\ref{pressure_four_diff_eq}, of the incident wave in Equation~\ref{four_diff_eq} as follows: 
\begin{equation}\label{shear_four_diff_eq}
\mathcal{F}_{1,2}\vecu(\xibo,\pm r_M,\omega)=\sqrt{\frac{\pi}{2}}\left[\hat{\vecG}_se^{i\kappa_sr_M}\mathcal{F}\vecf_s(\xibo,\pm\kappa_s-k_s,\omega)
                                                                            +\hat{\vecG}_pe^{i\kappa_pr_M}\mathcal{F}\vecf_s(\xibo,\pm\kappa_p-k_s,\omega)\right]                                                                                                                                                 
\end{equation}
\begin{equation}\label{pressure_four_diff_eq}
 \mathcal{F}_{1,2}\vecu(\xibo,\pm r_M,\omega)=\sqrt{\frac{\pi}{2}}\left[\hat{\vecG}_se^{i\kappa_sr_M}\mathcal{F}\vecf_p(\xibo,\pm\kappa_s-k_p,\omega)
                                                                           +\hat{\vecG}_pe^{i\kappa_pr_M}\mathcal{F}\vecf_p(\xibo,\pm\kappa_p-k_p,\omega)\right]
\end{equation}
Besides, the scattered wave (i.e. measured data) is a combination of both shear and pressure waves having different wavenumbers and amplitude vectors. The decomposition of the displacement measurements into its longitudinal and transversal parts is well studied in the literature, e.g. see the work by~\cite{DevOri86}. They used the properties of the S- and P-waves (see Remark~\ref{homo_elas_wave_remark}) and specific filters defined by the propagation vectors. In our case, one can simply use the fact that the vectors $\vecq_s$ and $\hat{\vecG}_s\vecg$ are orthogonal (i.e. $\vecq_s\cdot\hat{\vecG}_s\vecg=0$) and the vectors $\vecq_s$ and $\hat{\vecG}_s\vecg$ are parallel (i.e. $\vecq_p\times\hat{\vecG}_p\vecg=\bm{0}$), where the propagation vectors $\vecq_s$ and $\vecq_p$ defined by Equation~\ref{eq:prop_vect} and the partial Fourier transform of the Green's tensor $\hat{\vecG}_s$ and $\hat{\vecG}_p$ defined by Equation~\ref{eq:ft_green_tensor}. Finally, one gets the following formulae: 
\begin{itemize}
\item SS mode (S-image from S-excitation) is given by taking the curl product with the propagation vector $\vecq_p$ in Equation~\ref{shear_four_diff_eq}
\begin{equation}\label{ss_mode_eq}
 \vecq_p\times\mathcal{F}_{1,2}\vecu(\xibo,\pm r_M,\omega)=\sqrt{\frac{\pi}{2}}e^{i\kappa_sr_M}\vecq_p\times\hat{\vecG}_s\mathcal{F}\vecf_s(\xibo,\pm\kappa_s-k_s,\omega)
\end{equation}
\item PS mode (P-image from S-excitation) is given by taking the dot product with the propagation vector $\vecq_s$ in Equation~\ref{shear_four_diff_eq} 
\begin{equation}\label{ps_mode_eq}
 \vecq_s\cdot\mathcal{F}_{1,2}\vecu(\xibo,\pm r_M,\omega)=\sqrt{\frac{\pi}{2}}e^{i\kappa_pr_M}\vecq_s\cdot\hat{\vecG}_p\mathcal{F}\vecf_s(\xibo,\pm\kappa_p-k_s,\omega)
\end{equation}
\item SP mode (S-image from P-excitation) is given by taking the curl product with the propagation vector $\vecq_p$ in Equation~\ref{pressure_four_diff_eq} 
\begin{equation}\label{sp_mode_eq}
 \vecq_p\times\mathcal{F}_{1,2}\vecu(\xibo,\pm r_M,\omega)=\sqrt{\frac{\pi}{2}}e^{i\kappa_sr_M}\vecq_p\times\hat{\vecG}_s\mathcal{F}\vecf_p(\xibo,\pm\kappa_s-k_p,\omega)
\end{equation}
\item PP mode (P-image from P-excitation) is given by taking the dot product with the propagation vector $\vecq_s$ in Equation~\ref{pressure_four_diff_eq}
\begin{equation}\label{pp_mode_eq}
 \vecq_s\cdot\mathcal{F}_{1,2}\vecu(\xibo,\pm r_M,\omega)=\sqrt{\frac{\pi}{2}}e^{i\kappa_pr_M}\vecq_s\cdot\hat{\vecG}_p\mathcal{F}\vecf_p(\xibo,\pm\kappa_p-k_p,\omega)
\end{equation}
\end{itemize}
\begin{remark}
\label{back_prop_remark} 
There is an interference between different back-propagating waves at the medium discontinuities. For example, a back-propagated S-wave produces a S-wave as well as an induced P-wave (SP mode) in the back-propagation of S-wave and vice versa. Therefore, SP and PS modes are not the same due to different wavenumbers and amplitude vectors. 
\end{remark}
\begin{remark}
For an anisotropic medium, one needs to use the eigenvectors of Christoffel equation for the wave-field decomposition, for more details see~\cite{ZhaMcm10}.  
\end{remark}
\begin{remark}
\label{four_diff_remark}
Without rotation of the object, the measurements in both transmission and reflection imaging provide access to the shear scattering potential $\vecf_s$ on the two hemispheres in k-space with midpoint $-k_s\,\vece_3$ and radii $k_s$ and $k_p$, respectively
\begin{itemize}
\item SS mode
\[
S_{(0,0,-k_s),k_s}=\{(\xibo,\pm\kappa_s-k_s):\xibo\in\R^2,\,\xi_1^2+\xi_2^2<k_s^2\}
\]
\item PS mode 
\[
S_{(0,0,-k_s),k_p}=\{(\xibo,\pm\kappa_p-k_s):\xibo\in\R^2,\,\xi_1^2+\xi_2^2<k_p^2\}
\]
\end{itemize}
and the pressure scattering potential $\vecf_p$ on the two hemispheres in k-space with midpoint $-k_p\,\vece_3$ and radii $k_s$ and $k_p$, respectively
\begin{itemize}
\item SP mode 
\[
S_{(0,0,-k_p),k_s}=\{(\xibo,\pm\kappa_s-k_p):\xibo\in\R^2,\,\xi_1^2+\xi_2^2<k_s^2\}
\]
\item PP mode 
\[
S_{(0,0,-k_p),k_p}=\{(\xibo,\pm\kappa_p-k_p):\xibo\in\R^2,\,\xi_1^2+\xi_2^2<k_p^2\}
\]
\end{itemize}
One can conclude that $S_{(0,0,-k_s),k_p}\subset S_{(0,0,-k_s),k_s}$ and $S_{(0,0,-k_p),k_p}\subset S_{(0,0,-k_p),k_s}$ from the fact that $k_p<k_s$ (i.e. $v_p>v_s$, see Remark~\ref{homo_elas_wave_remark}). Then, it will be enough to measure the shear part of the displacement field to get the maximum coverage of the k-space. The different cases are depicted in Figure~\ref{four_diff_figb} for transmission imaging and in Figure~\ref{four_diff_figc} for reflection imaging.
\begin{figure}[h!]
\centering
\begin{tikzpicture}[scale=.8]
	\node[draw] at (-3,3) {PP mode};	
	\draw[->] (2,-3) -- (2,3) node[anchor = west]{{\color{black}$\xi_1$}};
	\draw[->] (3,1) -- (-1,-3) node[anchor = north]{{\color{black}$\xi_2$}};
	\draw[->] (-3.5,0) -- (3,0) node[anchor = north]{{\color{black}$\xi_3$}};
	\draw (0,-0.08) -- (0,0.08);
	\draw (0,-2) arc (-90:90:2cm);	
	\draw (0,2) arc (90:270:.5cm and 2cm);
	\draw[dashed] (0,-2) arc (-90:90:.5cm and 2cm);
	 \begin{scope}
	 \clip(3,3) -- (0,3) -- (0,2) arc (90:270:.5cm and 2cm) -- (0,-3) -- (3,-3) -- (3,3);
	 \shade[ball color=orange!60!white, opacity=0.70,overlay] (0,0) circle (2cm);
	 \end{scope}
	\node[anchor = south] at (0,0){$-k_p$};
\end{tikzpicture}
\hfil
\begin{tikzpicture}[scale=.8]	
	\node[draw] at (-3,3) {PS mode};
	\draw[->] (2,-3) -- (2,3) node[anchor = west]{{\color{black}$\xi_1$}};
	\draw[->] (3,1) -- (-1,-3) node[anchor = north]{{\color{black}$\xi_2$}};
	\draw[->] (-3.5,0) -- (3,0) node[anchor = north]{{\color{black}$\xi_3$}};
	\draw (-0.5,-0.08) -- (-0.5,0.08);
	\draw (-0.5,-2) arc (-90:90:2cm);	
	\draw (-0.5,2) arc (90:270:.5cm and 2cm);
	\draw[dashed] (-0.5,-2) arc (-90:90:.5cm and 2cm);
	\begin{scope}
	 \clip(1.5,3) -- (-0.5,3) -- (-0.5,2) arc (90:270:.5cm and 2cm) -- (-0.5,-3) -- (1.5,-3) -- (1.5,3);
	 \shade[ball color=orange!60!white, opacity=0.70] (-0.5,0) circle (2cm);
         \end{scope}
 	\node[anchor = south] at (-0.5,0){$-k_s$};
\end{tikzpicture}
\begin{tikzpicture}[scale=.8]	
	\node[draw] at (-3,3) {SP mode};	
	\draw[->] (2,-3) -- (2,3) node[anchor = west]{{\color{black}$\xi_1$}};
	\draw[->] (3,1) -- (-1,-3) node[anchor = north]{{\color{black}$\xi_2$}};
	\draw[->] (-3.5,0) -- (3,0) node[anchor = north]{{\color{black}$\xi_3$}};
	\draw (0,-0.08) -- (0,0.08);
	\draw (0,-2.5) arc (-90:90:2.5cm);	
	\draw (0,2.5) arc (90:270:.5cm and 2.5cm);
	\draw[dashed] (0,-2.5) arc (-90:90:.5cm and 2.5cm);
	 \begin{scope}
	 \clip(3.5,3.5) -- (0,3.5) -- (0,2.5) arc (90:270:.5cm and 2.5cm) -- (0,-3.5) -- (3.5,-3.5) -- (3.5,3.5);
	 \shade[ball color=orange!60!white, opacity=0.70,overlay] (0,0) circle (2.5cm);
	 \end{scope}
	\node[anchor = south] at (0,0){$-k_p$};
\end{tikzpicture}
\hfil
\begin{tikzpicture}[scale=.8]	
	\node[draw] at (-3,3) {SS mode};	
	\draw[->] (2,-3) -- (2,3) node[anchor = west]{{\color{black}$\xi_1$}};
	\draw[->] (3,1) -- (-1,-3) node[anchor = north]{{\color{black}$\xi_2$}};
	\draw[->] (-3.5,0) -- (3,0) node[anchor = north]{{\color{black}$\xi_3$}};
	\draw (-0.5,-0.08) -- (-0.5,0.08);
	\draw (-0.5,-2.5) arc (-90:90:2.5cm);	
	\draw (-0.5,2.5) arc (90:270:.5cm and 2.5cm);
	\draw[dashed] (-0.5,-2.5) arc (-90:90:.5cm and 2.5cm);
	\begin{scope}
	 \clip(3,3.5) -- (-0.5,3.5) -- (-0.5,2.5) arc (90:270:.5cm and 2.5cm) -- (-0.5,-3.5) -- (3,-3.5) -- (3,3.5);
	 \shade[ball color=orange!60!white, opacity=0.70,overlay] (-0.5,0) circle (2.5cm);
	 \end{scope}
 	\node[anchor = south] at (-0.5,0){$-k_s$};
\end{tikzpicture}
\caption{Accessible points in k-space for transmission imaging.}
\label{four_diff_figb}
\end{figure}
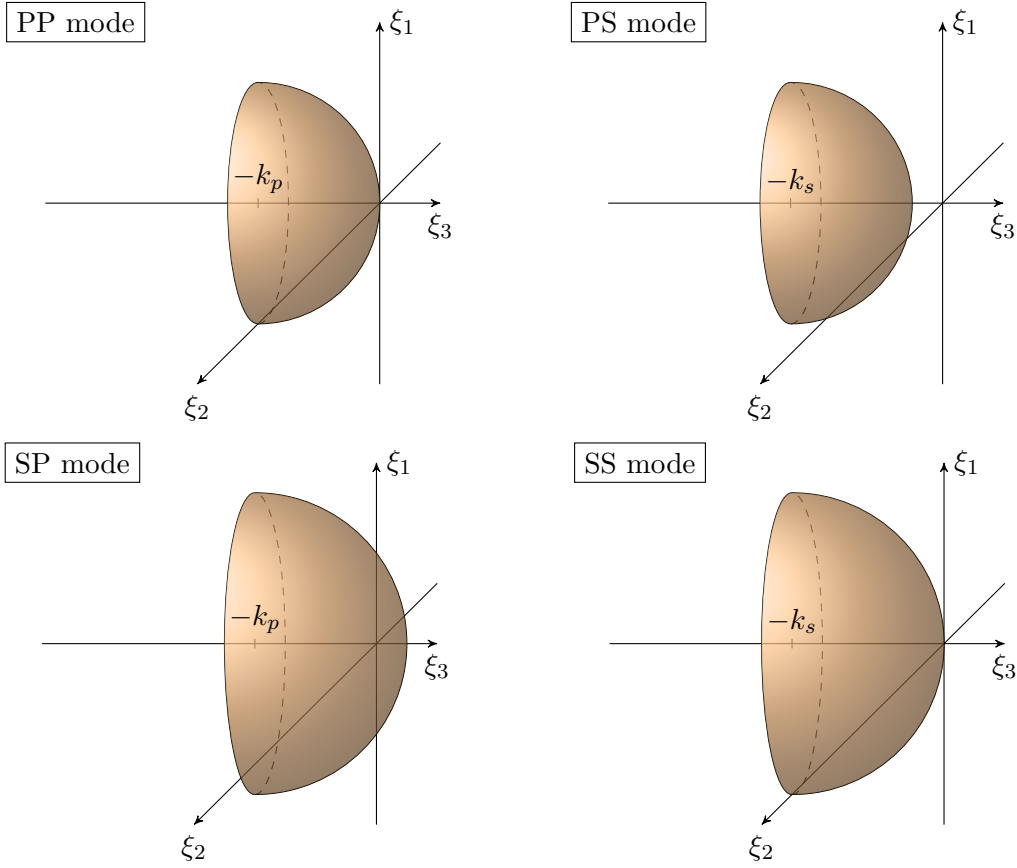
\begin{figure}[h!]
\centering
\begin{tikzpicture}[scale=.8]
	\node[draw] at (-3,3) {PP mode};
	\draw (-2.5,0) -- (-.5,0);
	\draw (-3.5,0) -- (-2.5,0);
	\draw (0,2) arc (90:270:2cm);	
	\draw (0,2) arc (90:270:.5cm and 2cm);
	\draw (0,-2) arc (-90:90:.5cm and 2cm);
	\begin{scope}
	\clip(-2.5,3) -- (0,3) -- (0,2) arc (90:270:.5cm and 2cm) -- (0,-3) -- (-2.5,-3) -- (-2.5,3);
	\shade[ball color=orange!60!white, opacity=0.70] (0,0) circle (2cm);
	\end{scope}
	\shade[top color = orange!60!black, bottom color = orange!20!white, opacity=0.70] (.5,0) arc (0:360:.5cm and 2cm);
	\draw[->] (2,-3) -- (2,3) node[anchor = west]{{\color{black}$\xi_1$}}; 	
 	\draw[->] (3,1) -- (-1,-3) node[anchor = north]{{\color{black}$\xi_2$}};
	\draw[->] (-.5,0) -- (3,0) node[anchor = north]{{\color{black}$\xi_3$}};
	\draw (0,-0.08) -- (0,0.08);
 	\node[anchor = south] at (0,0){$-k_p$};
\end{tikzpicture}
\hfil
\begin{tikzpicture}[scale=.8]
        \node[draw] at (-3,3) {PS mode};
	\draw (-2.5,0) -- (-1.,0);
	\draw (-0.5,2) arc (90:270:2cm);	
	\draw (-0.5,2) arc (90:270:.5cm and 2cm);
	\draw (-0.5,-2) arc (-90:90:.5cm and 2cm);
	\begin{scope}
	\clip(-2.5,3) -- (-0.5,3) -- (-0.5,2) arc (90:270:.5cm and 2cm) -- (-0.5,-3) -- (-2.5,-3) -- (-2.5,3);
	\shade[ball color=orange!60!white, opacity=0.70] (-0.5,0) circle (2cm);
	\end{scope}
	\shade[top color = orange!60!black, bottom color = orange!20!white, opacity=0.70] (0,0) arc (0:360:.5cm and 2cm);
	\draw[->] (2,-3) -- (2,3) node[anchor = west]{{\color{black}$\xi_1$}}; 	
 	\draw[->] (3,1) -- (-1,-3) node[anchor = north]{{\color{black}$\xi_2$}};
	\draw[->] (-1.,0) -- (3,0) node[anchor = north]{{\color{black}$\xi_3$}};
	\draw (-3.5,0) -- (-2.5,0);
	\draw (-0.5,-0.08) -- (-0.5,0.08);
	\node[anchor = south] at (-0.5,0){$-k_s$};
\end{tikzpicture}	
\begin{tikzpicture}[scale=.8]	
	\node[draw] at (-3,3) {SP mode};
	\draw (-3,0) -- (-0.5,0);
	\draw (0,2.5) arc (90:270:2.5cm);	
	\draw (0,2.5) arc (90:270:.5cm and 2.5cm);
	\draw (0,-2.5) arc (-90:90:.5cm and 2.5cm);
	\begin{scope}
	\clip(-2.5,3.5) -- (0,3.5) -- (0,2.5) arc (90:270:.5cm and 2.5cm) -- (0,-3.5) -- (-2.5,-3.5) -- (-2.5,3.5);
	\shade[ball color=orange!60!white, opacity=0.70] (0,0) circle (2.5cm);
	\end{scope}
	\shade[top color = orange!60!black, bottom color = orange!20!white, opacity=0.70] (.5,0) arc (0:360:.5cm and 2.5cm);
	\draw[->] (2,-3) -- (2,3) node[anchor = west]{{\color{black}$\xi_1$}}; 	
 	\draw[->] (3,1) -- (-1,-3) node[anchor = north]{{\color{black}$\xi_2$}};
	\draw[->] (-.5,0) -- (3,0) node[anchor = north]{{\color{black}$\xi_3$}};
	\draw (-3.5,0) -- (-3,0);
	\draw (0,-0.08) -- (0,0.08);
 	\node[anchor = south] at (0,0){$-k_p$};
\end{tikzpicture}
\hfil
\begin{tikzpicture}[scale=.8]
	\node[draw] at (-3,3) {SS mode};
	\draw (-3,0) -- (-1,0);
	\draw (-0.5,2.5) arc (90:270:2.5cm);	
	\draw (-0.5,2.5) arc (90:270:.5cm and 2.5cm);
	\draw (-0.5,-2.5) arc (-90:90:.5cm and 2.5cm);
	\begin{scope}
	\clip(-3,3.5) -- (-0.5,3.5) -- (-0.5,2.5) arc (90:270:.5cm and 2.5cm) -- (-0.5,-3.5) -- (-3,-3.5) -- (-3,3.5);
	\shade[ball color=orange!60!white, opacity=0.70] (-0.5,0) circle (2.5cm);
	\end{scope}
	\shade[top color = orange!60!black, bottom color = orange!20!white, opacity=0.70] (0,0) arc (0:360:.5cm and 2.5cm);
	\draw[->] (2,-3) -- (2,3) node[anchor = west]{{\color{black}$\xi_1$}}; 	
 	\draw[->] (3,1) -- (-1,-3) node[anchor = north]{{\color{black}$\xi_2$}};
	\draw[->] (-1,0) -- (3,0) node[anchor = north]{{\color{black}$\xi_3$}};
	\draw (-3.5,0) -- (-3,0);
	\draw (-0.5,-0.08) -- (-0.5,0.08);
	\node[anchor = south] at (-0.5,0){$-k_s$};
\end{tikzpicture}		
\caption{Accessible points in k-space for reflection imaging.}
\label{four_diff_figc}
\end{figure}
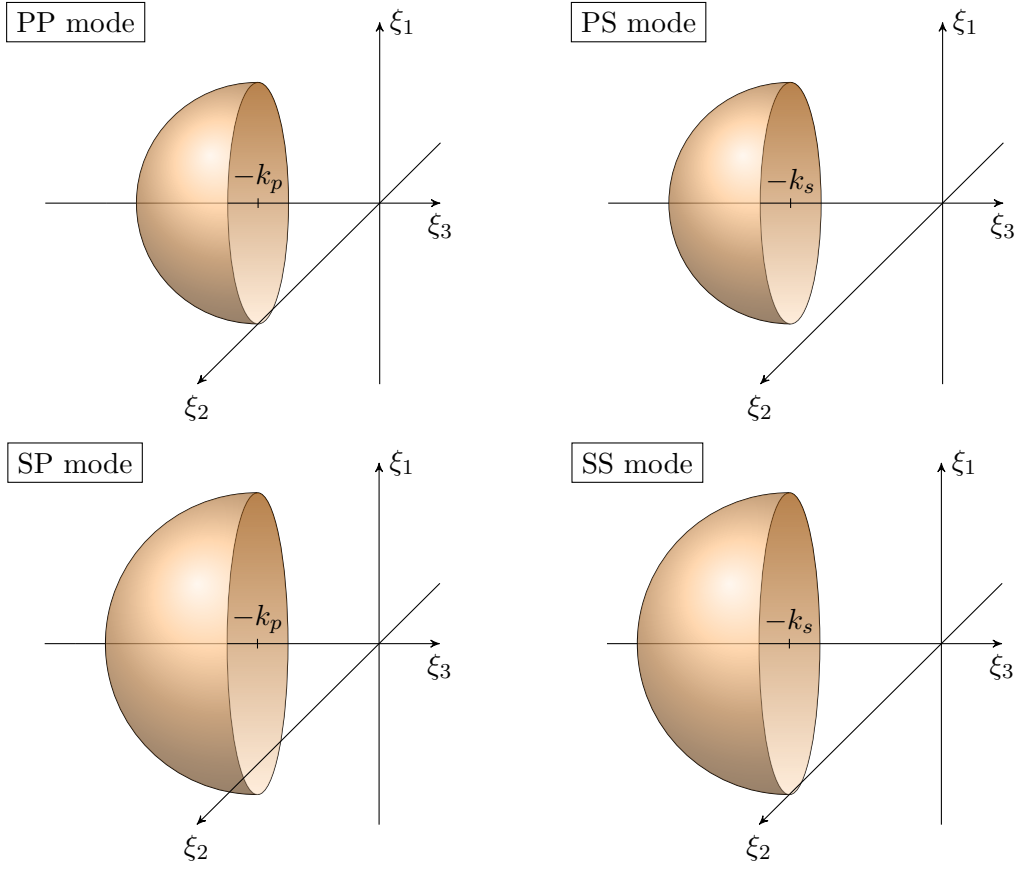
\end{remark}
\begin{remark}
Blackledge \textit{et al.}~\cite{BlaBurHopWom87} presents the Fourier slice diffraction theorem where the Fourier transform of the projection yields the Fourier transform of the object along a plane (i.e. 2D measure). However, the Fourier diffraction corollary~\ref{four_diff_coro} relates the 2D Fourier transform of the measured scattered field projection with the Fourier transform of the object along a hemisphere (i.e. 3D measure) in the spatial frequency domain. 
\end{remark}

\section{Wave field inversion}
\label{wave_field_inversion}
The next step is the inversion process for the different modes. However, one remarks that the matrices $\hat{\vecG}_s$ and $\hat{\vecG}_p$ are singular, Equation~\ref{eq:ft_green_tensor}, so we can't directly invert the formulae~\ref{ss_mode_eq}, \ref{ps_mode_eq}, \ref{sp_mode_eq} and \ref{pp_mode_eq}. Then, one needs to develop the right-hand side of each formula by taking into account the properties of the shear/pressure amplitude vectors $(\mathtt{A}_s,\mathtt{A}_p)$, the propagation vectors $(\vecq_s,\vecq_p)$ and the partial Fourier transform of Green's tensor $(\hat{\vecG}_s,\hat{\vecG}_p)$. 

Firstly, one presents the Fourier transform of the shear scattering potential $\vecf_s$, Equation~\ref{fs_eq}, as follows 
\begin{equation}\label{eq:ft_fs}
\mathcal{F}\vecf_s(\xibo,\pm\kappa_\alpha-k_s,\omega)=\left[k_s^2\hat{\delta\mu}-ik_s(\pm\kappa_\alpha-k_s)\hat{\delta\mu}-\omega^2\hat{\delta\rho}\right]\mathtt{A}_s-ik_s\hat{\delta\mu}(\xibo'\cdot\mathtt{A}_s)\vece_3
\end{equation}
and the Fourier transform of the pressure scattering potential $\vecf_p$, Equation~\ref{fp_eq}, as follows 
\begin{equation}\label{eq:ft_fp}
\mathcal{F}\vecf_p(\xibo,\pm\kappa_\alpha-k_p,\omega)=\left[k_p^2(\hat{\delta\lambda}+2\hat{\delta\mu})-ik_p(\pm\kappa_\alpha-k_p)(\hat{\delta\lambda}+2\hat{\delta\mu})-\omega^2\hat{\delta\rho}\right]\mathtt{A}_p-ik_p\hat{\delta\lambda}(\vece_3\cdot\mathtt{A}_p)\xibo'
\end{equation}
where the parameters $\hat{\delta\mu}$, $\hat{\delta\lambda}$ and $\hat{\delta\rho}$ are, respectively, the Fourier transform of the perturbed elastic parameters $\delta\mu$, $\delta\lambda$ and $\delta\rho$ defined at $(\xibo,\pm\kappa_\alpha-k_\beta)$ in the the Fourier domain, for $\alpha,\beta=s,p$. 
\begin{lemma} PP mode: Let the assumptions of Corollary~\ref{four_diff_coro} be satisfied. One defines the PP scattering function $\hat{f}_{pp_{(\hat{\delta\mu},\hat{\delta\lambda},\hat{\delta\rho})}}$ in the Schwartz space $\mathcal{S}(\R^3)$ by 
\[
\hat{f}_{pp}(\xibo,\pm\kappa_p-k_p,\omega):=k_p^2\kappa_p(\hat{\delta\lambda}+2\hat{\delta\mu})-ik_p\kappa_p(\pm\kappa_p-k_p)(\hat{\delta\lambda}+2\hat{\delta\mu})-\omega^2\kappa_p\hat{\delta\rho}-k_p\|\xibo'\|\hat{\delta\lambda}
\]
Then, one gets 
\begin{equation}\label{eq:ppmode}
\hat{f}_{pp}(\xibo,\pm\kappa_p-k_p,\omega)=\sqrt{\frac{2}{\pi}}e^{-i\kappa_pr_M}\frac{\vecq_s\cdot\mathcal{F}_{1,2}\vecu(\xibo,\pm r_M,\omega)}{(\vece_3\cdot\mathtt{A}_p)\vecq_s\cdot\vecq_p}
\end{equation}
\end{lemma}
\begin{proof}
We derive the form $\hat{\vecG}_p\mathcal{F}\vecf_p$, using Equation~\ref{eq:ft_fp}, $\hat{\vecG}_p\mathtt{A}_p=\kappa_p(\vece_3\cdot\mathtt{A}_p)\vecq_p$ and $\hat{\vecG}_p\xibo'=-i\|\xibo'\|\vecq_p$. One gets 
\[ 
\hat{\vecG}_p\mathcal{F}\vecf_p(\xibo,\pm\kappa_p-k_p,\omega)=\underbrace{\left[k_p^2\kappa_p(\hat{\delta\lambda}+2\hat{\delta\mu})-ik_p\kappa_p(\pm\kappa_p-k_p)(\hat{\delta\lambda}+2\hat{\delta\mu})-\omega^2\kappa_p\hat{\delta\rho}-k_p\|\xibo'\|\hat{\delta\lambda}\right]}_{\hat{f}_{pp}(\xibo,\pm\kappa_p-k_p,\omega)}(\vece_3\cdot\mathtt{A}_p)\vecq_p
\]
Then, PP mode, Equation~\ref{pp_mode_eq} has the following formula  
\[
 \vecq_s\cdot\mathcal{F}_{1,2}\vecu(\xibo,\pm r_M,\omega)=\sqrt{\frac{\pi}{2}}e^{i\kappa_pr_M}\hat{f}_{pp}(\xibo,\pm\kappa_p-k_p,\omega)(\vece_3\cdot\mathtt{A}_p)\vecq_s\cdot\vecq_p
\]
\end{proof}
\begin{lemma} PS mode: Let the assumptions of Corollary~\ref{four_diff_coro} be satisfied. One defines the PS scattering function $\hat{f}_{ps_{(\hat{\delta\mu},\hat{\delta\rho})}}$ in the Schwartz space $\mathcal{S}(\R^3)$ by 
\[
\hat{f}_{ps}(\xibo,\pm\kappa_p-k_s,\omega):=-i\left[k_s^2\hat{\delta\mu}-ik_s(\pm\kappa_p-k_s)\hat{\delta\mu}-\omega^2\hat{\delta\rho}+k_s\kappa_p\hat{\delta\mu}\right]
\]
Then, one gets 
\begin{equation}\label{eq:psmode}
\hat{f}_{ps}(\xibo,\pm\kappa_p-k_s,\omega)=\sqrt{\frac{2}{\pi}}e^{-i\kappa_pr_M}\frac{\vecq_s\cdot\mathcal{F}_{1,2}\vecu(\xibo,\pm r_M,\omega)}{(\xibo'\cdot\mathtt{A}_s)\vecq_s\cdot\vecq_p}
\end{equation}
\end{lemma}
\begin{proof}
Here, we represent the term $\hat{\vecG}_p\mathcal{F}\vecf_s$, with the fact that $\hat{\vecG}_p\mathtt{A}_s=-i(\xibo'\cdot\mathtt{A}_s)\vecq_p$ and $\hat{\vecG}_p\vece_3=\kappa_p\vecq_p$, in Equation~\ref{eq:ft_fs}, as follows 
\[
\hat{\vecG}_p\mathcal{F}\vecf_s(\xibo,\pm\kappa_p-k_s,\omega)=\underbrace{-i\left[k_s^2\hat{\delta\mu}-ik_s(\pm\kappa_\alpha-k_s)\hat{\delta\mu}-\omega^2\hat{\delta\rho}+k_s\kappa_p\hat{\delta\mu}\right]}_{\hat{f}_{ps}(\xibo,\pm\kappa_p-k_s,\omega)}(\xibo'\cdot\mathtt{A}_s)\vecq_p
\]
Likewise, PS mode, Equation~\ref{ps_mode_eq} has the following formula  
\[
 \vecq_s\cdot\mathcal{F}_{1,2}\vecu(\xibo,\pm r_M,\omega)=\sqrt{\frac{\pi}{2}}e^{i\kappa_pr_M}\hat{f}_{ps}(\xibo,\pm\kappa_p-k_s,\omega)(\xibo'\cdot\mathtt{A}_s)\vecq_s\cdot\vecq_p
\]
\end{proof}
\begin{lemma} SP mode: Let the assumptions of Corollary~\ref{four_diff_coro} be satisfied. One defines the SP scattering function $\hat{f}_{sp_{(\hat{\delta\mu},\hat{\delta\lambda},\hat{\delta\rho})}}$ in the Schwartz space $\mathcal{S}(\R^3)$ by 
\[
\hat{f}_{sp}(\xibo,\pm\kappa_s-k_p,\omega):=(\kappa_s^2-\kappa_s\kappa_p+k_s^2)\left[(k_p^2-ik_p(\pm\kappa_s-k_p))(\hat{\delta\lambda}+2\hat{\delta\mu})-\omega^2\hat{\delta\rho}-k_p\kappa_s\hat{\delta\lambda}\right]
\]
Then, one gets two formula for SH horizontal shear and SV vertical shear-image from P-wave incidence as follows 
\begin{equation}\label{eq:spmode}
\hat{f}_{sp}(\xibo,\pm\kappa_s-k_p,\omega)=\sqrt{\frac{2}{\pi}}e^{-i\kappa_sr_M}\frac{\vecq_p\times\mathcal{F}_{1,2}\vecu(\xibo,\pm r_M,\omega)\cdot\vece_i}{(\vece_3\cdot\mathtt{A}_p)\vece_3\times\vecq_p\cdot\vece_i},\quad\text{for }\;i=1,2
\end{equation}
\end{lemma}
\begin{proof}
Now, we establish the term $\hat{\vecG}_s\mathcal{F}\vecf_p$, with the fact that $\hat{\vecG}_s\mathtt{A}_p=\kappa_s(\vece_3\cdot\mathtt{A}_p)\vecq_s+k_s^2\mathtt{A}_p$ and $\hat{\vecG}_s\xibo'=-i\kappa_s^2\vecq_s-i\kappa_sk_s^2\vece_3$, in Equation~\ref{eq:ft_fp}, as follows 
\begin{equation*}
\begin{split}
\hat{\vecG}_s\mathcal{F}\vecf_p(\xibo,\pm\kappa_s-k_p,\omega)
  &=\left[k_p^2\kappa_s(\hat{\delta\lambda}+2\hat{\delta\mu})-ik_p\kappa_s(\pm\kappa_s-k_p)(\hat{\delta\lambda}+2\hat{\delta\mu})-\omega^2\kappa_s\hat{\delta\rho}-k_p\kappa_s^2\hat{\delta\lambda}\right](\vece_3\cdot\mathtt{A}_p)\vecq_s
\\&+\left[k_p^2k_s^2(\hat{\delta\lambda}+2\hat{\delta\mu})-ik_pk_s^2(\pm\kappa_s-k_p)(\hat{\delta\lambda}+2\hat{\delta\mu})-\omega^2k_s^2\hat{\delta\rho}-k_pk_s^2\kappa_s\hat{\delta\lambda}\right]\mathtt{A}_p
\end{split}
\end{equation*}
Next, we derive the term $\vecq_p\times\hat{\vecG}_s\mathcal{F}\vecf_p$, with the fact that $\vecq_p\times\vecq_s=(\kappa_s-\kappa_p)\vece_3\times\vecq_p$ and $\vecq_p\times\mathtt{A}_p=(\vece_3\cdot\mathtt{A}_p)\vece_3\times\vecq_p$, in the last equation, as follows 
\[
\hat{\vecG}_s\mathcal{F}\vecf_p(\xibo,\pm\kappa_s-k_p,\omega)=\underbrace{(\kappa_s^2-\kappa_s\kappa_p+k_s^2)\left[(k_p^2-ik_p(\pm\kappa_s-k_p))(\hat{\delta\lambda}+2\hat{\delta\mu})-\omega^2\hat{\delta\rho}-k_p\kappa_s\hat{\delta\lambda}\right]}_{\hat{f}_{sp}(\xibo,\pm\kappa_s-k_p,\omega)}(\vece_3\cdot\mathtt{A}_p)\vece_3\times\vecq_p
\]
Afterwards, SP mode, Equation~\ref{sp_mode_eq} has the following form 
\[
 \vecq_p\times\mathcal{F}_{1,2}\vecu(\xibo,\pm r_M,\omega)=\sqrt{\frac{\pi}{2}}e^{i\kappa_sr_M}\hat{f}_{sp}(\xibo,\pm\kappa_s-k_p,\omega)(\vece_3\cdot\mathtt{A}_p)\vece_3\times\vecq_p
\]
\end{proof}
\begin{lemma} SS mode: Let the assumptions of Corollary~\ref{four_diff_coro} be satisfied. One defines the first SS scattering function $\hat{f}_{ss,1_{(\hat{\delta\mu},\hat{\delta\rho})}}$ in the Schwartz space $\mathcal{S}(\R^3)$ by 
\[
\hat{f}_{ss,1}(\xibo,\pm\kappa_s-k_s,\omega):=-i(\kappa_s-\kappa_p)\left[k_s^2\hat{\delta\mu}-ik_s(\pm\kappa_s-k_s)\hat{\delta\mu}-\omega^2\hat{\delta\rho}+k_s\kappa_s\hat{\delta\mu}\right]-ik_s^3\hat{\delta\mu}
\]
Then, one has 
\begin{equation}\label{eq:ss1mode}
\hat{f}_{ss,1}(\xibo,\pm\kappa_s-k_s,\omega)=\sqrt{\frac{2}{\pi}}\frac{e^{-i\kappa_sr_M}}{d'}(\vecq_p\times\mathcal{F}_{1,2}\vecu(\xibo,\pm r_M,\omega))\cdot\vece_1(\vece_3\times\mathtt{A}_s)\cdot\vece_2
\end{equation}
and the second SS scattering function $\hat{f}_{ss,2_{(\hat{\delta\mu},\hat{\delta\rho})}}$ in the Schwartz space $\mathcal{S}(\R^3)$ defined by 
 \[
\hat{f}_{ss,2}(\xibo,\pm\kappa_s-k_s,\omega):=-k_s^2\kappa_s\left[k_s^2\hat{\delta\mu}-ik_s(\pm\kappa_s-k_s)\hat{\delta\mu}-\omega^2\hat{\delta\rho}\right]
\]
Then, one gets 
\begin{equation}\label{eq:ss2mode}
\hat{f}_{ss,2}(\xibo,\pm\kappa_s-k_s,\omega)=\sqrt{\frac{2}{\pi}}\frac{e^{-i\kappa_sr_M}d''}{d'}(\vecq_p\times\mathcal{F}_{1,2}\vecu(\xibo,\pm r_M,\omega))\cdot\vece_1
\end{equation}
where $d''= (\sqrt{\frac{2}{\pi}}e^{-i\kappa_sr_M}(\vecq_p\times\mathcal{F}_{1,2}\vecu(\xibo,\pm r_M,\omega))\cdot\vece_2-(\xibo'\cdot\mathtt{A}_s)(\vece_3\times\vecq_p)\cdot\vece_2)$
and $d'=(\xibo'\cdot\mathtt{A}_s)(\vece_3\times\vecq_p)\cdot\vece_1(\vece_3\times\mathtt{A}_s)\cdot\vece_2+(\vece_3\times\mathtt{A}_s)\cdot\vece_1d''$.
\end{lemma}
\begin{proof}
We represent the form $\hat{\vecG}_s\mathcal{F}\vecf_s$, with the fact that $\hat{\vecG}_s\mathtt{A}_s=-i(\xibo'\cdot\mathtt{A}_s)\vecq_s+k_s^2\mathtt{A}_s$ and $\hat{\vecG}_s\vece_3=\kappa_s\vecq_s+k_s^2\vece_3$, in Equation~\ref{eq:ft_fs}, as follows 
\begin{equation*}
\begin{split}
\hat{\vecG}_s\mathcal{F}\vecf_s(\xibo,\pm\kappa_s-k_s,\omega)
  &=-i\left[k_s^2\hat{\delta\mu}-ik_s(\pm\kappa_s-k_s)\hat{\delta\mu}-\omega^2\hat{\delta\rho}+k_s\kappa_s\hat{\delta\mu}\right](\xibo'\cdot\mathtt{A}_s)\vecq_s-ik_s^3\hat{\delta\mu}(\xibo'\cdot\mathtt{A}_s)\vece_3
\\&+k_s^2\left[k_s^2\hat{\delta\mu}-ik_s(\pm\kappa_s-k_s)\hat{\delta\mu}-\omega^2\hat{\delta\rho}\right]\mathtt{A}_s
\end{split}
\end{equation*}
Next, we establish the form $\vecq_p\times\hat{\vecG}_s\mathcal{F}\vecf_s$, with the fact that $\vecq_p\times\vecq_s=(\kappa_s-\kappa_p)\vece_3\times\vecq_p$ and $\vecq_p\times\mathtt{A}_s=-\kappa_s\vece_3\times\mathtt{A}_s$, in the last equation, as follows 
\begin{equation*}
\begin{split}
\vecq_p\times\hat{\vecG}_s\mathcal{F}\vecf_s(\xibo,\pm\kappa_s-k_s,\omega)&=\underbrace{-k_s^2\kappa_s\left[k_s^2\hat{\delta\mu}-ik_s(\pm\kappa_s-k_s)\hat{\delta\mu}-\omega^2\hat{\delta\rho}\right]}_{\hat{f}_{ss,1}(\xibo,\pm\kappa_s-k_s,\omega)}\vece_3\times\mathtt{A}_s
  \\&\underbrace{-\left(i(\kappa_s-\kappa_p)\left[k_s^2\hat{\delta\mu}-ik_s(\pm\kappa_s-k_s)\hat{\delta\mu}-\omega^2\hat{\delta\rho}+k_s\kappa_s\hat{\delta\mu}\right]+ik_s^3\hat{\delta\mu}\right)}_{\hat{f}_{ss,2}(\xibo,\pm\kappa_s-k_s,\omega)}(\xibo'\cdot\mathtt{A}_s)\vece_3\times\vecq_p
\end{split}
\end{equation*}
Solving the algebraic form, one gets 
\[
\hat{f}_{ss,1}d=(\vecq_p\times\hat{\vecG}_s\mathcal{F}\vecf_s)\cdot\vece_1(\vece_3\times\mathtt{A}_s)\cdot\vece_2
\]
and
\[
\hat{f}_{ss,2}d=(\vecq_p\times\hat{\vecG}_s\mathcal{F}\vecf_s)\cdot\vece_1((\vecq_p\times\hat{\vecG}_s\mathcal{F}\vecf_s)\cdot\vece_2-(\xibo'\cdot\mathtt{A}_s)(\vece_3\times\vecq_p)\cdot\vece_2)
\]
where
\[
d=(\xibo'\cdot\mathtt{A}_s)(\vece_3\times\vecq_p)\cdot\vece_1(\vece_3\times\mathtt{A}_s)\cdot\vece_2+(\vece_3\times\mathtt{A}_s)\cdot\vece_1((\vecq_p\times\hat{\vecG}_s\mathcal{F}\vecf_s)\cdot\vece_2-(\xibo'\cdot\mathtt{A}_s)(\vece_3\times\vecq_p)\cdot\vece_2)
\]
Finally, SS mode, Equation~\ref{ss_mode_eq} has the following forms
\[
\hat{f}_{ss,1}(\xibo,\pm\kappa_s-k_s,\omega)d'=\sqrt{\frac{2}{\pi}}e^{-i\kappa_sr_M}(\vecq_p\times\mathcal{F}_{1,2}\vecu(\xibo,\pm r_M,\omega))\cdot\vece_1(\vece_3\times\mathtt{A}_s)\cdot\vece_2
\]
and
\begin{equation*}
\begin{split}
\hat{f}_{ss,2}(\xibo,\pm\kappa_s-k_s,\omega)
d'&=\sqrt{\frac{2}{\pi}}e^{-i\kappa_sr_M}(\vecq_p\times\mathcal{F}_{1,2}\vecu(\xibo,\pm r_M,\omega))\cdot\vece_1(\sqrt{\frac{2}{\pi}}e^{-i\kappa_sr_M}
\\&(\vecq_p\times\mathcal{F}_{1,2}\vecu(\xibo,\pm r_M,\omega))\cdot\vece_2-(\xibo'\cdot\mathtt{A}_s)(\vece_3\times\vecq_p)\cdot\vece_2)
\end{split}
\end{equation*}
\end{proof}
\begin{remark}
\begin{itemize}
\item In the previous literature~\cite{BlaBurHopWom87,BeyBur90}, one needs to measure the scattered wave-field decomposition (i.e. $\vecu_{pp}$, $\vecu_{ps}$, $\vecu_{sp}$ and $\vecu_{ss}$) to reconstruct the elastic parameters for each mode. Experimentally, It is a challenging task. However, in our case one needs to measure only the full scattered displacement field for the multimodes.       
\item In the previous literature~\cite{BlaBurHopWom87,BeyBur90}, the authors used at most three position of source/receiver to reconstruct the multi-parameters (i.e. Lamé parameters, mass density). In this work, we will use the information from at most three different modes. Another alternative is to use measurements from transmission and reflection acquisitions.
\end{itemize}
\end{remark}

\section{K-space coverage}
\label{kspace_coverage}
One can extract k-space informations from the scattered displacement field using a single direction according to the Fourier diffraction Theorem~\ref{four_diff_theo}, as depicted in Figure~\ref{four_diff_figb} and Figure~\ref{four_diff_figc}. Actually, different coverage of k-space can be obtained in terms of angular diversity by varying the illumination direction for a fixed frequency or in terms of frequency diversity by imposing an impulsive illumination for a frequency band. Up to our knowledge, the k-space coverage wasn't treated in the previous literature for the elastic case unlike the acoustic case studied in~\cite{Langenberg1987,KirQueRitSchSet21,FauKirQueSchSet21_report} striking similarity to the PP mode. The aim of this section is to present various band limited versions of frequency coverage.   
\subsection{Angular diversity}
\label{angular_diversity}
Now, we will present the back-propagation formula for a rotated observation space. To this end, one gives the matrix of a rotation $R$ by an angle $\theta=\theta(t):[0,L]\to\R$ around the axis $\vecn=\vecn(t):[0,L]\to\S^2$ as follows 
\[
    R_{\theta,\vecn}:=\left(
    \begin{array}{ccc} 
    n_1^2(1-c)+c            & n_1n_2(1-c)-n_3s    & n_1n_3(1-c)+n_2s    \\
    n_1n_2(1-c)+n_3s    & n_2^2(1-c)+c           & n_2n_3(1-c)-n_1s     \\
    n_1n_3(1-c)-n_2s     & n_2n_3(1-c)+n_1s   & n_3^2(1-c)+c            \\
    \end{array}
    \right)
\]
where $\vecn=(n_1,n_2,n_3)$ a unit vector, $c:=\cos{\theta}$ and $s:=\sin{\theta}$. Assuming that the scattering object undergoes a continuous rotation $R$ with varying rotation axis. 

Then, one denotes by $\vecutt,\;0\leqslant t\leqslant L,$ the scattered displacement field of the rotated object. Under Born's approximation, it satisfies 
\[
    \left\{\begin{aligned}
    & \sigmabott(\x,\omega)=\vecC^0(\x):\varepsilonbo\left[\vecutt\right](\x,\omega)                      \\
    & \nablabo\cdot\sigmabott(\x,\omega)+\omega^2\rho^0(\x)\vecutt(\x,\omega) = \vecf_s(\x,\omega)\circ R_{\theta(t),\vecn(t)}e^{ik_sx_3} + \vecf_p(\x,\omega)\circ R_{\theta(t),\vecn(t)}e^{ik_px_3}             
    \end{aligned} \right.
\]
The full set of measurements in the transmission and reflection setup, respectively, is then given by 
\[
\vecutt(x_1,x_2,\pm r_M,\omega),\;\;x_1,\,x_2\in\R,\;0\leqslant t\leqslant L, 
\]
and, according to formulae~\ref{eq:ppmode},~\ref{eq:psmode},~\ref{eq:spmode},~\ref{eq:ss1mode} and,~\ref{eq:ss2mode}, it is related to the scattering potentials $f_{pp}$, $f_{ps}$, $f_{sp}$, $f_{ss,1}$ and $f_{ss,2}$, respectively, defined in the rotated observation space. 

Denoting by 
\[
\mathcal{U}_{\alpha}:=\{(\xibo,t)\in\R^3,\,\xi_1^2+\xi_2^2<k_\alpha^2,\;0\leqslant t\leqslant L\},\;\;\text{for }\alpha=s,p
\]
the sets where Equation~\ref{eq:ppmode}, Equation~\ref{eq:psmode}, Equation~\ref{eq:spmode}, Equation~\ref{eq:ss1mode} and Equation~\ref{eq:ss2mode} are valid and can be used for the reconstruction purpose, recall Remark~\ref{four_diff_remark}. Moreover, the map that traces out the accessible domains in k-space is denoted by 
\[
T^\pm_{\alpha\beta}:\mathcal{U}_{\alpha}\to\R^3,\;\;T^\pm_{\alpha\beta}(\xibo,t):=R_{\theta(t),\vecn(t)}(\xibo,\pm\kappa_\alpha-k_\beta)^\top,\;\;\text{for }\alpha,\beta=s,p
\]
In the reconstruction formula below, we have to take into account the number of times a point $\y$ in k-space is covered by $T^\pm_{\alpha\beta}$. This number, sometimes referred to as {\it Banach indicatrix} of $T^\pm_{\alpha\beta}$, will be denoted by Card$((T^\pm_{\alpha\beta})^{-1})$, where Card(A) is the cardinality of a set A. Finally, the approximation of $f$ we wish to reconstruct is 
\[
f^\pm_{\text{bp}}(\vecr):=(2\pi)^{-\frac32}\int_{T^\pm_{\alpha\beta}(\mathcal{U}_\alpha)}e^{i\y\cdot\vecr}\mathcal{F}f(\y)\dy.
\]
The set $T^\pm_{\alpha\beta}(\mathcal{U}_\alpha)$ will be referred to as the {\it frequency coverage} or {\it k-space coverage} of the experimental setup for multimodes.
\begin{theorem}\label{back_prop_theorem}
Let the assumptions of Corollary~\ref{four_diff_coro} be satisfied. In addition, assume that $\theta\in C^1[0,L]$ and $\vecn\in(C^1[0,L],\S^2)$. Then, for all $\vecr\in\R^3$, one has
\begin{itemize}
\item PP mode 
\[
f^{\pm}_{\text{pp}}(\vecr,\omega)=\frac{1}{2\pi^2}\int_{\mathcal{U}_p}e^{iT^\pm_{pp}(\xibo,t)\cdot\vecr}\frac{\vecq_s\cdot\mathcal{F}_{1,2}\vecutt(\xibo,\pm r_M,\omega)}{(\vece_3\cdot\mathtt{A}_p)\vecq_s\cdot\vecq_p}\frac{e^{-i\kappa_pr_M}|\nabla T^\pm_{pp}(\xibo,t)|}{\text{Card}(T^\pm_{pp}(T^\pm_{pp}(\xibo,t))^{-1})}\,\text{d}(\xibo,t)
\]
\item PS Mode 
\[
f^\pm_{\text{ps}}(\vecr,\omega)=\frac{1}{2\pi^2}\int_{\mathcal{U}_p}e^{iT^\pm_{ps}(\xibo,t)\cdot\vecr}\frac{\vecq_s\cdot\mathcal{F}_{1,2}\vecutt(\xibo,\pm r_M,\omega)}{(\xibo'\cdot\mathtt{A}_s)\vecq_s\cdot\vecq_p}
\frac{e^{-i\kappa_pr_M}|\nabla T^\pm_{ps}(\xibo,t)|}{\text{Card}(T^\pm_{ps}(T^\pm_{ps}(\xibo,t))^{-1})}\,\text{d}(\xibo,t)
\]
\item SP mode 
\[
f^\pm_{\text{sp}}(\vecr,\omega)=\frac{1}{2\pi^2}\int_{\mathcal{U}_s}e^{iT^\pm_{sp}(\xibo,t)\cdot\vecr}\frac{\vecq_p\times\mathcal{F}_{1,2}\vecutt(\xibo,\pm r_M,\omega)\cdot\vece_i}{(\vece_3\cdot\mathtt{A}_p)\vece_3\times\vecq_p\cdot\vece_i}\frac{e^{-i\kappa_sr_M}|\nabla T^\pm_{sp}(\xibo,t)|}{\text{Card}(T^\pm_{sp}(T^\pm_{sp}(\xibo,t))^{-1})}\,\text{d}(\xibo,t),\;\text{for }i=1,2
\]
\item SS mode 
\[
f^\pm_{\text{ss,1}}(\vecr,\omega)=\frac{1}{2\pi^2}\int_{\mathcal{U}_s}e^{iT^\pm_{ss}(\xibo,t)\cdot\vecr}\frac{(\vecq_p\times\mathcal{F}_{1,2}\vecutt(\xibo,\pm r_M,\omega))\cdot\vece_1(\vece_3\times\mathtt{A}_s)\cdot\vece_2}{d'}\frac{e^{-i\kappa_sr_M}|\nabla T^\pm_{ss}(\xibo,t)|}{\text{Card}(T^\pm_{ss}(T^\pm_{ss}(\xibo,t))^{-1})}\,\text{d}(\xibo,t),
\]
\[
f^\pm_{\text{ss,2}}(\vecr,\omega)=\frac{1}{2\pi^2}\int_{\mathcal{U}_s}e^{iT^\pm_{ss}(\xibo,t)\cdot\vecr}
\frac{(\vecq_p\times\mathcal{F}_{1,2}\vecutt(\xibo,\pm r_M,\omega))\cdot\vece_1d''}{d'}\frac{e^{-i\kappa_sr_M}|\nabla T^\pm_{ss}(\xibo,t)|}{\text{Card}(T^\pm_{ss}(T^\pm_{ss}(\xibo,t))^{-1})}\,\text{d}(\xibo,t),
\]
\end{itemize}
where the scattering potentials $f_{pp}$, $f_{ps}$, $f_{sp}$, $f_{ss,1}$ and $f_{ss,2}$ are defined in Section~\ref{wave_field_inversion} and $|\nabla T^\pm_{\alpha\beta}|$ is the magnitude of the Jacobian determinant of $T^\pm_{\alpha\beta}$, for $\alpha,\beta=s,p$.
\end{theorem}
\begin{proof}
The proof has the same lines as the one presented in the paper by~\cite{KirQueRitSchSet21} except the inversion step developed in Section~\ref{wave_field_inversion}.
\end{proof}
Recalling the following corollary.  
\begin {corollary}
\label{back_prop_coro}
~\cite{KirQueRitSchSet21}
Let the assumptions of Theorem~\ref{back_prop_theorem} be satisfied. If $\vecn'(t)=\bm{0}$, then 
\[
|\nabla T^\pm_{\alpha\beta}(\xibo,t|=\frac{k_\beta|\theta'(t)||n_2\xi_1-n_1\xi_2|}{\kappa_\alpha}
\]
If in addition $\vecn\neq\vece_3$ and $\theta$ is strictly increasing with $\theta(0)=0$ and $\theta(L)=2\pi$, then 
\[
\text{Card}(T^\pm_{\alpha\beta}(T^\pm_{\alpha\beta}(\xibo,t))^{-1})=2
\]
for almost every $(\xibo,t)\in\mathcal{U}_\alpha$.
\end{corollary}
\begin{remark}
The back-propagation formulae depends on the k-space cardinality, which is delicate to calculate. It doesn't change from the acoustic to the elastic cases, as it is always presented by a hemisphere with different midpoints and radii. 
\end{remark}
\begin{example} 
Full uniform rotation around the $x_1$-axis: Consider a rotation around the $x_1$-axis with the rotation matrix 
\[
    R_{\theta,\vece_1}:=\left(
    \begin{array}{ccc} 
    1           & 0    & 0    \\
    0    & \cos{\theta}           & -\sin{\theta}     \\
    0     & \sin{\theta}   & \cos{\theta}            \\
    \end{array}
    \right)
\]
and $\theta(t)=t$, for $t\in[0,2\pi]$. From Corollary~\ref{back_prop_coro}, one gets
\[
f^\pm_{\alpha\beta}(\vecr,\omega)=\frac{k_\beta}{4\pi^2}\int_{\mathcal{U}_\alpha}\frac{1}{\kappa_\alpha}e^{i(R_{\theta,\vece_1}\vech^\pm_{\alpha\beta}\cdot\vecr-\kappa_\alpha r_M)}\hat{f}_{\alpha\beta}(\xibo,\pm\kappa_\alpha-k_\beta,\omega)|\xi_2|\,\text{d}(\xibo,\theta)
\]
where $\vech^\pm_{\alpha\beta}:=(\xibo,\pm\kappa_\alpha-k_\beta)^\top$, for $\alpha,\beta=s,p$.
\label{examp1}
\end{example}
Our next step is to illustrate the corresponding sets $R_{\theta,\vece_1}\vech^\pm_{\alpha\beta}$. Consider $\y\in T^+_{\alpha\beta}(\mathcal{U}_\alpha)$, i.e.
 $\y=(y_1,y_2,y_3)^\top=R_{\theta,\vece_1}(\xi_1,\xi_2,\kappa_\alpha-k_\beta)^\top$, where $\xi_1^2+\xi_2^2<k_\alpha^2$ and $0\leq\theta\leq2\pi$. Then, it holds that $y_1=\xi_1$ and 
\[
y_2^2+y_3^2=k_\alpha^2+k_\beta^2-2k_\beta\sqrt{k_\alpha^2-\xi_1^2-\xi_2^2}-\xi_1^2
\]
Since $k_\beta>0$, one has $y_2^2+y_3^2<k_\alpha^2+k_\beta^2-\xi_1^2$, in other form $\|\y\|^2<k_\alpha^2+k_\beta^2$. Since $\xi_2^2\in[0,k_\alpha^2-\xi_1^2)$, this implies 
\[
y_2^2+y_3^2\geq k_\alpha^2+k_\beta^2-2k_\beta\sqrt{k_\alpha^2-\xi_1^2}-\xi_1^2=\left(k_\beta-\sqrt{k_\alpha^2-\xi_1^2}\right)^2
\]
with some rearranging gives $y_1^2+\left(\|(y_2,y_3)\|-k_\beta\right)^2\leq k_\alpha^2$. In summary, one gets 
\[
T^+_{\alpha\beta}(\mathcal{U}_\alpha)=\{\y\in\R^3:\|\y\|^2<k_\alpha^2+k_\beta^2,\,y_1^2+\left(\|(y_2,y_3)\|-k_\beta\right)^2\leq k_\alpha^2\},\;\text{for }\alpha,\beta=s,p
\]
\begin{itemize}
\item PP mode: K-space is covered by a horn torus of radius $k_p$ given by 
\[
T^+_{pp}(\mathcal{U}_p)=\{\y\in\R^3:\|\y\|^2<2k_p^2,\,y_1^2+\left(\|(y_2,y_3)\|-k_p\right)^2\leq k_p^2\}
\]                                                    
\item PS mode: K-space is covered by a ring torus, with $k_p$ is the radius of the tube and $k_s$ is the radius from the center of the hole to the center of the tours tube, given by 
\[
T^+_{ps}(\mathcal{U}_p)=\{\y\in\R^3:\|\y\|^2<k_p^2+k_s^2,\,y_1^2+\left(\|(y_2,y_3)\|-k_s\right)^2\leq k_p^2\}
\]
\item SP mode: K-space is covered by a spindle torus, with $k_s$ is the radius of the tube and $k_p$ is the radius from the center of the hole to the center of the tours tube, given by 
\[
T^+_{sp}(\mathcal{U}_s)=\{\y\in\R^3:\|\y\|^2<k_s^2+k_p^2,\,y_1^2+\left(\|(y_2,y_3)\|-k_p\right)^2\leq k_s^2\}
\]
\item SS mode: K-space is covered by a horn torus of radius $k_s$ given by 
\[
T^+_{ss}(\mathcal{U}_s)=\{\y\in\R^3:\|\y\|^2<2k_s^2,\,y_1^2+\left(\|(y_2,y_3)\|-k_s\right)^2\leq k_s^2\}
\]
\end{itemize}
The different cases are depicted in Figure~\ref{fig:angdiver} for a full uniform rotation about the $x_1$-axis. 
\begin{figure}[h!]
\centering
\begin{tikzpicture}[scale=0.6]
	\node[draw] at (-3.5,5) {PP mode};
   	\draw[->]   (0,-4.5) -- (0,4.5) node[right,black] {$y_1$};
	\draw (0,0) circle (2);
    	\draw (0,0) circle (1.41);
    	\draw (1,0) circle (1);
   	\draw (-1,0) circle (1);
    	\fill[cyan!70!blue,opacity=.2] (1,1) arc (90:270:1) -- (1,-1) arc (-45:45:1.41);
	\fill[cyan!70!blue,opacity=.2] (-1,-1) arc (-90:90:1) -- (-1,1) arc (135:225:1.41);
   	\fill[magenta,opacity=.2] (1,1) arc (45:-45:1.41) -- (1,-1) arc (-90:90:1);
   	\fill[magenta,opacity=.2] (-1,-1) arc (225:135:1.41) -- (-1,1) arc (90:270:1);
	\draw[->] (0,0) -- (-.7,1.217);
	\node[below] at (-.7,2) {$\sqrt{2}k_p$};
	\draw[->] (1,0) -- (1.5,-.866);
	\node[below] at (1.5,-0.9) {$k_p$};
	\draw[->] (0,0) -- (1,1.732);
	\node[below] at (1.2,2.5) {$2k_p$};
	\node[above] at (.75,0) {TI};
	\node[above] at (-.75,0) {TI};
	\node[below] at (1.6,0) {RI};
	\node[below] at (-1.6,0) {RI};
\end{tikzpicture}
\hfil
\begin{tikzpicture}[scale=0.6]
	\node[draw] at (-3.5,5) {PS mode};
    	\draw[->] (0,-4.5) -- (0,4.5) node[right,black] {$y_1$};
	\draw (0,0) circle (3); 
    	\draw (0,0) circle (2.23);
    	\draw (1.75,0) circle (1.25);
    	\draw (-1.75,0) circle (1.25);
	\fill[cyan!70!blue,opacity=.2] (1.849,1.246) arc (85:275:1.25) -- (1.849,-1.246) arc (-34:34:2.23); 
	\fill[cyan!70!blue,opacity=.2] (-1.849,1.246) arc (95:-95:1.25) -- (-1.849,1.246) arc (146:214:2.23);
   	\fill[magenta,opacity=.2] (1.849,1.246) arc (34:-34:2.23) -- (1.849,-1.246) arc (-85:85:1.25);
     	\fill[magenta,opacity=.2] (-1.849,1.246) arc (146:214:2.23) -- (-1.849,-1.246) arc (265:95:1.25);
   	\draw[->] (0,0) -- (-1.03,1.978);
   	\node[below] at (-1.03,3) {$\sqrt{k_s^2+k_p^2}$};
   	\draw[->] (1.75,0) -- (2.5,-1);
   	\node[below] at (2.5,-1.2) {$\frac{3k_p+k_s}{4}$};
   	\draw[->] (0,0) -- (1,2.83);
   	\node[below] at (1.2,3.7) {$k_s+k_p$};
	\node[above] at (1.5,0) {TI};
	\node[above] at (-1.5,0) {TI};
	\node[below] at (2.6,0) {RI};
	\node[below] at (-2.6,0) {RI};
\end{tikzpicture}

\begin{tikzpicture}[scale=0.6]
	\node[draw] at (-3.5,5) {SP mode};
  	\fill[cyan!70!blue,opacity=.2] (1.389,1.744) arc (85:275:1.75) -- (1.389,-1.744) arc (-51:45:2.23);
   	\fill[cyan!70!blue,opacity=.2] (-1.389,-1.744) arc (-95:95:1.75) -- (-1.389,1.744) arc (129:225:2.23);
        \fill[white] (0,1.22) arc (45:-45:1.73) ;
      	\fill[white] (0,-1.22) arc (225:135:1.73) ;
   	\fill[magenta,opacity=.2] (1.389,1.744) arc (50:-50:2.28) -- (1.389,-1.744) arc (-85:85:1.75);
    	\fill[magenta,opacity=.2] (-1.389,-1.744) arc (230:130:2.28) -- (-1.389,1.744) arc (95:260:1.75);
    	\draw[->] (0,-4.5) -- (0,4.5) node[right,black] {$y_1$};
	\draw (0,0) circle (3); 
    	\draw (0,0) circle (2.23);
    	\draw (1.25,0) circle (1.75);
    	\draw (-1.25,0) circle (1.75);
   	\draw[->] (0,0) -- (-1.03,1.978);
   	\node[below] at (-1.03,3) {$\sqrt{k_s^2+k_p^2}$};
   	\draw[->] (1.25,0) -- (2,-1.581);
   	\node[below] at (2,-1.6) {$\frac{3k_s+k_p}{4}$};
   	\draw[->]   (0,0) -- (1,2.83);
   	\node[below] at (1.2,3.7) {$k_s+k_p$};
	\node[above] at (1.5,0) {TI};
	\node[above] at (-1.5,0) {TI};
	\node[below] at (2.6,0) {RI};
	\node[below] at (-2.6,0) {RI};
\end{tikzpicture}
\hfil
\begin{tikzpicture}[scale=0.6]
	\node[draw] at (-3.5,5) {SS mode};
	\draw[->] (0,-4.5) -- (0,4.5) node[right,black] {$y_1$};
	\draw (0,0) circle (4);
    	\draw (0,0) circle (2.83);
    	\draw (2,0) circle (2);
    	\draw (-2,0) circle (2);
       	\fill[cyan!70!blue,opacity=.2] (2,2) arc (90:270:2) -- (2,-2) arc (-45:45:2.83);
        \fill[cyan!70!blue,opacity=.2] (-2,-2) arc (-90:90:2) -- (-2,2) arc (135:225:2.83);
   	\fill[magenta,opacity=.2] (2,2) arc (45:-45:2.83) -- (2,-2) arc (-90:90:2);
    	\fill[magenta,opacity=.2] (-2,-2) arc (225:135:2.83) -- (-2,2) arc (90:270:2);
   	\draw[->] (0,0) -- (-1,2.644);
    	\node[below] at (-1,3.5) {$\sqrt{2}k_s$};
   	\draw[->]   (2,0) -- (3,-1.732);
    	\node[below] at (3,-1.8) {$k_s$};
   	\draw[->] (0,0) -- (1,3.873);
   	\node[below] at (1.2,4.6) {$2k_s$};
	\node[above] at (1.5,0) {TI};
	\node[above] at (-1.5,0) {TI};
	\node[below] at (3.5,0) {RI};
	\node[below] at (-3.5,0) {RI};
\end{tikzpicture}
\caption{K-space coverage for a full uniform rotation about the $x_1$-axis: cross-section of transmission imaging (TI) (magenta part)- cross-section of reflection imaging (RI) (cyan part).}  
\label{fig:angdiver}
\end{figure}
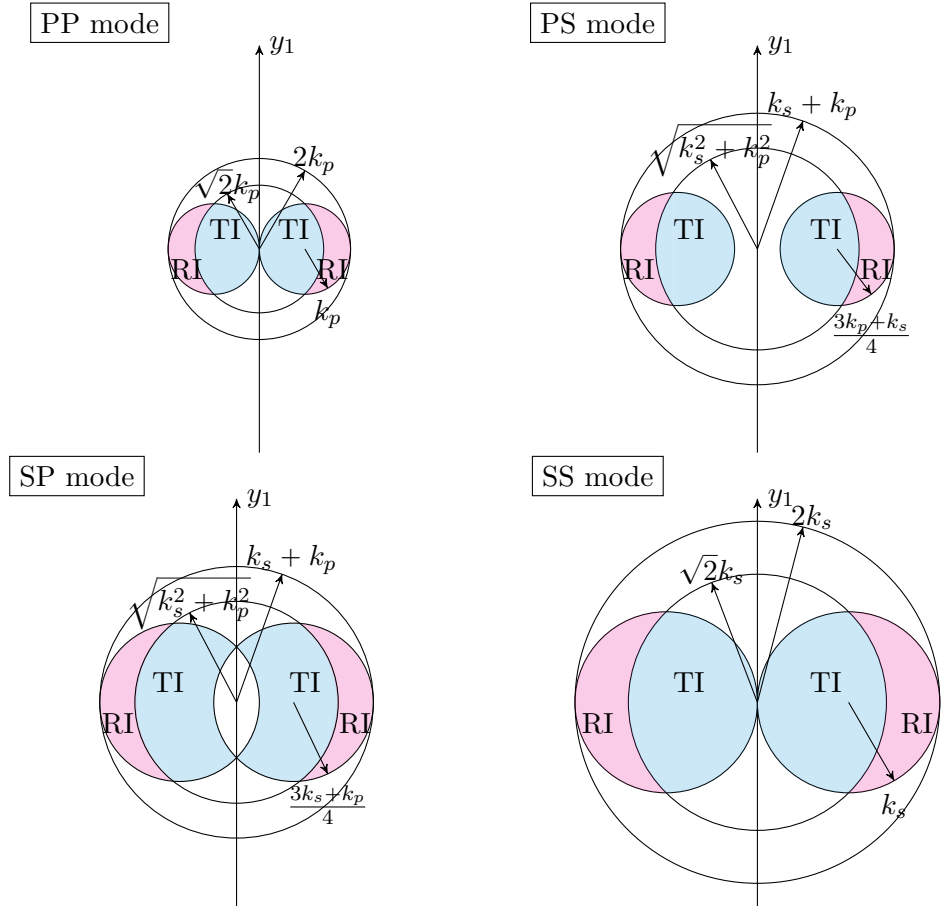\begin{remark}
One remarks that the transmission imaging offers better spatial resolution for each modes by dint of the larger k-space coverage. 
\end{remark}

\subsection{Frequency diversity}
\label{frequency_diversity} 
Here, one considers a mono-angle multifrequencies plane wave excitation by varying the shear wavenumber $k_s$ over a set $K_s=\left[k_{s,\mathrm{min}},k_{s,\mathrm{max}}\right]$ and the pressure wavenumber $k_p$ over a set $K_p=\left[k_{p,\mathrm{min}},k_{p,\mathrm{max}}\right]$.  

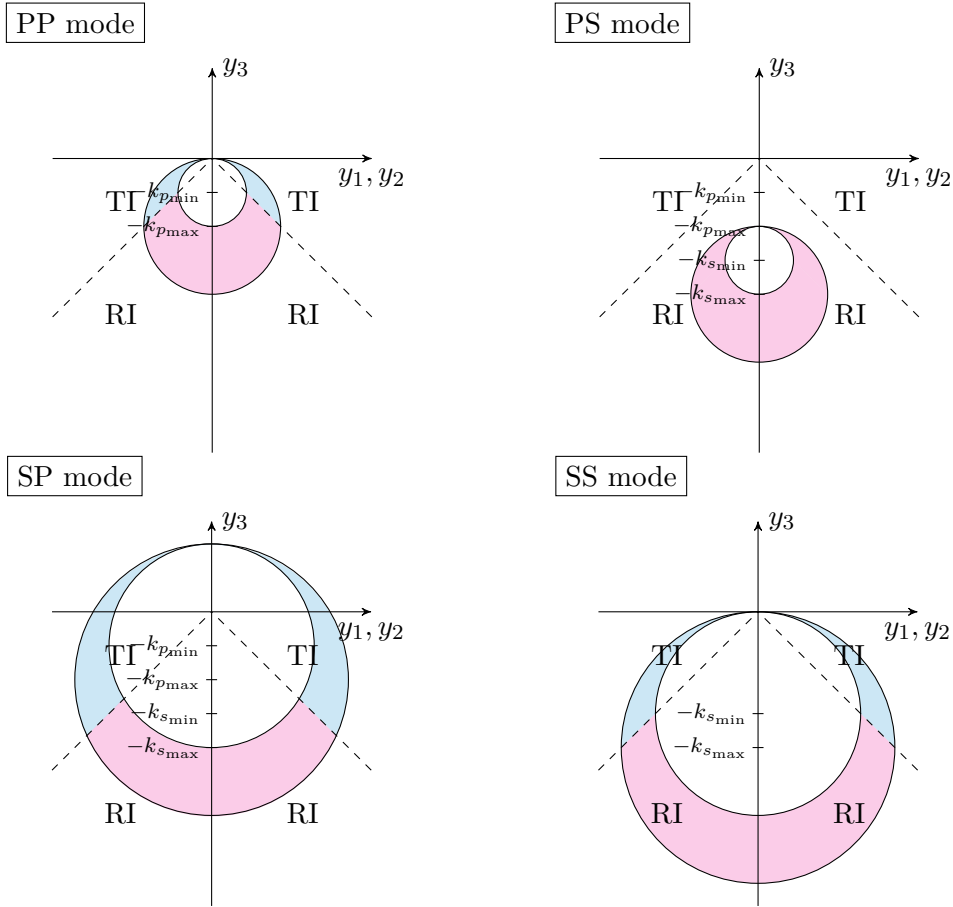
\begin{figure}[h!]
\centering
\begin{tikzpicture}[scale=0.6]
	\node[draw] at (-3,3) {PP mode};
	\fill[cyan!70!blue,opacity=.2] (-1.5,-1.5) arc (180:0:1.5) -- (.75,-.75) arc (0:180:.75);
	\fill[magenta,opacity=.2] (-1.5,-1.5) arc (180:360:1.5) -- (.75,-.75) arc (0:-180:.75);
	
	\draw[black] (0,-.75) circle (.75cm);
	\draw[black] (0,-1.5) circle (1.5cm);
	
	\draw[->]   (-3.5,0) -- (3.5,0) node[below,black] {$y_1,y_2$};
	\draw[->]   (0,-6.5) -- (0,2) node[right,black] {$y_3$};
	
	\draw[black,dashed]   (-3.5,-3.5) -- (0,0);
	\draw[black,dashed]   (3.5,-3.5) -- (0,0);
	\node[below] at (2,-0.5) {TI};
	\node[below] at (-2,-0.5) {TI};
	\node[below] at (2,-3) {RI};
	\node[below] at (-2,-3) {RI};
	
	\draw (-0.12,-.75) -- (0.12,-.75);
	\node[left] at (0,-0.75) {\scriptsize{$-k_{p_{\mathrm{min}}}$}};
	\draw (-0.12,-1.5) -- (0.12,-1.5);
	\node[left] at (0,-1.5) {\scriptsize{$-k_{p_{\mathrm{max}}}$}};
\end{tikzpicture}
\hfil
\begin{tikzpicture}[scale=0.6]
	\node[draw] at (-3,3) {PS mode};
	\fill[magenta,opacity=.2] (-1.5,-3) arc (180:0:1.5) -- (.75,-2.25) arc (0:180:.75);
	\fill[magenta,opacity=.2] (-1.5,-3) arc (180:360:1.5) -- (.75,-2.25) arc (0:-180:.75);
	
	\draw[black] (0,-2.25) circle (.75cm);
	\draw[black] (0,-3) circle (1.5cm);
	
	\draw[->]   (-3.5,0) -- (3.5,0) node[below,black] {$y_1,y_2$};
	\draw[->]   (0,-6.5) -- (0,2) node[right,black] {$y_3$};
	
	\draw[black,dashed]   (-3.5,-3.5) -- (0,0);
	\draw[black,dashed]   (3.5,-3.5) -- (0,0);
	\node[below] at (2,-0.5) {TI};
	\node[below] at (-2,-0.5) {TI};
	\node[below] at (2,-3) {RI};
	\node[below] at (-2,-3) {RI};
	
	\draw (-0.12,-.75) -- (0.12,-.75);
	\node[left] at (0,-0.75) {\scriptsize{$-k_{p_{\mathrm{min}}}$}};
	\draw (-0.12,-1.5) -- (0.12,-1.5);
	\node[left] at (0,-1.5) {\scriptsize{$-k_{p_{\mathrm{max}}}$}};
	\draw[black] (-0.12,-2.25) -- (0.12,-2.25);
	\node[left] at (0,-2.25) {\scriptsize{$-k_{s_{\mathrm{min}}}$}};
	\draw[black] (-0.12,-3) -- (0.12,-3);
	\node[left] at (0,-3) {\scriptsize{$-k_{s_{\mathrm{max}}}$}};
\end{tikzpicture}

\begin{tikzpicture}[scale=0.6]
	\node[draw] at (-3,3) {SP mode};
	\fill[cyan!70!blue,opacity=.2] (-2.743,-2.743) arc (204:-24:3) -- (1.921,-1.921) arc (-31.24:211.24:2.25);
	\fill[magenta,opacity=.2] (-2.743,-2.743) arc (204:336:3) -- (1.921,-1.921) arc (329:211:2.25);
	
	\draw[black] (0,-.75) circle (2.25cm);
	\draw[black] (0,-1.5) circle (3cm);
	
	\draw[->]   (-3.5,0) -- (3.5,0) node[below,black] {$y_1,y_2$};
	\draw[->]   (0,-6.5) -- (0,2) node[right,black] {$y_3$};
	\draw[black,dashed]   (-3.5,-3.5) -- (0,0);
	\draw[black,dashed]   (3.5,-3.5) -- (0,0);
	\node[below] at (2,-0.5) {TI};
	\node[below] at (-2,-0.5) {TI};
	\node[below] at (2,-4) {RI};
	\node[below] at (-2,-4) {RI};
		
	\draw (-0.12,-.75) -- (0.12,-.75);
	\node[left] at (0,-0.75) {\scriptsize{$-k_{p_{\mathrm{min}}}$}};
	\draw (-0.12,-1.5) -- (0.12,-1.5);
	\node[left] at (0,-1.5) {\scriptsize{$-k_{p_{\mathrm{max}}}$}};
	\draw[black] (-0.12,-2.25) -- (0.12,-2.25);
	\node[left] at (0,-2.25) {\scriptsize{$-k_{s_{\mathrm{min}}}$}};
	\draw[black] (-0.12,-3) -- (0.12,-3);
	\node[left] at (0,-3) {\scriptsize{$-k_{s_{\mathrm{max}}}$}};
\end{tikzpicture}
\hfil
\begin{tikzpicture}[scale=0.6]
	\node[draw] at (-3,3) {SS mode};
	\fill[cyan!70!blue,opacity=.2] (-3,-3) arc (180:0:3) -- (2.25,-2.25) arc (0:180:2.25);
	\fill[magenta,opacity=.2] (-3,-3) arc (180:360:3) -- (2.25,-2.25) arc (0:-180:2.25);
	
	\draw[black] (0,-2.25) circle (2.25cm);
	\draw[black] (0,-3) circle (3cm);
	
	\draw[->]   (-3.5,0) -- (3.5,0) node[below,black] {$y_1,y_2$};
	\draw[->]   (0,-6.5) -- (0,2) node[right,black] {$y_3$};
	\draw[black,dashed]   (-3.5,-3.5) -- (0,0);
	\draw[black,dashed]   (3.5,-3.5) -- (0,0);
	\node[below] at (2,-0.5) {TI};
	\node[below] at (-2,-0.5) {TI};
	\node[below] at (2,-4) {RI};
	\node[below] at (-2,-4) {RI};		
	\draw (-0.12,-2.25) -- (0.12,-2.25);
	\node[left] at (0,-2.25) {\scriptsize{$-k_{s_{\mathrm{min}}}$}};
	\draw (-0.12,-3) -- (0.12,-3);
	\node[left] at (0,-3) {\scriptsize{$-k_{s_{\mathrm{max}}}$}};
\end{tikzpicture}
\caption{
K-space coverage where the shear wavenumber takes values in $[k_{s_{\mathrm{min}}},k_{s_{\mathrm{max}}}]$ and the pressure wavenumber takes values in $[k_{p_{\mathrm{min}}},k_{p_{\mathrm{max}}}]$: transmission imaging (TI) (magenta part) - reflection imaging (RI) (cyan part).} 
\label{fig:freqdiver}
\end{figure}

Denoting by $\vecu$ the scattered displacement field depending on the angular frequency $\omega$, for $\omega_{\mathrm{min}}\leqslant\omega\leqslant\omega_{\mathrm{max}}$. Under Born's approximation, it satisfies 
\[
	\left\{\begin{aligned}
    		& \sigmabo(\x,\omega)=\vecC^0(\x):\varepsilonbo\left[\vecu\right](\x,\omega)                      \\
    		& \nablabo\cdot\sigmabo(\x,\omega)+\omega^2\rho^0(\x)\vecu(\x,\omega) = \vecf_s(\x,\omega)e^{ik_sx_3} + 			\vecf_p(\x,\omega)e^{ik_px_3}             
    	\end{aligned} \right.
\]
The full set of measurements in the transmission and reflection setup, respectively, is then given by 
\[
	\vecu(x_1,x_2,\pm r_M,\omega),\;\;x_1,\,x_2\in\R,\;\;\omega_{\mathrm{min}}\leqslant\omega\leqslant\omega_{\mathrm{max}} 
\]
Notice that $\kappa_{\alpha}$ also depends on $\omega$. The resulting k-space coverage
\[
	\mathcal{Y}_{\alpha\beta} = \left\{(\xibo,\pm\kappa_\alpha-k_\beta)^\top\in\R^3: k_{\alpha,\beta}\in K_{\alpha,\beta}, \, \|\xibo\|< k_\alpha\right\},\;\;\text{for }\alpha,\beta=s,p
\]
is covered by the region between two internally tangent circles $(-k_{\alpha_\mathrm{max}}\vece_3,k_{\beta_\mathrm{max}})$ and $(-k_{\beta_\mathrm{min}}\vece_3,k_{\alpha_\mathrm{min}})$. Then, the k-space coverage consists of all points $\y\in\R^3$ such that $\|(y_1,y_2)\| \le k_{\alpha_\mathrm{max}}$ 
and
\[
	\sqrt{k_{\alpha_\mathrm{max}}^2-y_1^2-y_2^2} - k_{\alpha_\mathrm{max}}
	\ge y_3
	\ge 
	\begin{cases} 
		-\|(y_1,y_2)\|                                                                   &\|(y_1,y_2)\|\ge k_{\alpha_\mathrm{min}}\\
		\sqrt{k_{\alpha_\mathrm{min}}^2-y_1^2-y_2^2} - k_{\alpha_\mathrm{min}} & \text{otherwise}
	\end{cases}
\]
for $\alpha,\beta=s,p$. 
\begin{itemize}
\item PP mode: K-space is covered by the region between two internally tangent circles $(-k_{p_\mathrm{max}}\vece_3,k_{p_\mathrm{max}})$ and $(-k_{p_\mathrm{min}}\vece_3,k_{p_\mathrm{min}})$. It consists of all points $\y\in\R^3$ such that $\|(y_1,y_2)\| \le k_{p_\mathrm{max}}$ and
\[
	\sqrt{k_{p_\mathrm{max}}^2-y_1^2-y_2^2} - k_{p_\mathrm{max}}
	\ge y_3
	\ge 
	\begin{cases} 
		-\|(y_1,y_2)\|                                                                   &\|(y_1,y_2)\|\ge k_{p_\mathrm{min}}\\
		\sqrt{k_{p_\mathrm{min}}^2-y_1^2-y_2^2} - k_{p_\mathrm{min}} & \text{otherwise}
	\end{cases}
\] 
\item PS mode: K-space is covered by the region between two internally tangent circles $(-k_{s_\mathrm{max}}\vece_3,k_{p_\mathrm{max}})$ and $(-k_{s_\mathrm{min}}\vece_3,k_{p_\mathrm{min}})$. It consists of all points $\y\in\R^3$ such that $\|(y_1,y_2)\| \le k_{p_\mathrm{max}}$ and
\[
	\sqrt{k_{p_\mathrm{max}}^2-y_1^2-y_2^2} - k_{s_\mathrm{max}}
	\ge y_3
	\ge 
	\begin{cases} 
		-\|(y_1,y_2)\|                                                                   &\|(y_1,y_2)\|\ge k_{p_\mathrm{min}}\\
		\sqrt{k_{p_\mathrm{min}}^2-y_1^2-y_2^2} - k_{s_\mathrm{min}} & \text{otherwise}
	\end{cases}
\]
\item SP mode: K-space is covered by the region between two internally tangent circles $(-k_{p_\mathrm{max}}\vece_3,k_{s_\mathrm{max}})$ and $(-k_{p_\mathrm{min}}\vece_3,k_{s_\mathrm{min}})$. It consists of all points $\y\in\R^3$ such that $\|(y_1,y_2)\| \le k_{s_\mathrm{max}}$ and
\[
	\sqrt{k_{s_\mathrm{max}}^2-y_1^2-y_2^2} - k_{p_\mathrm{max}}
	\ge y_3
	\ge 
	\begin{cases} 
		-\|(y_1,y_2)\|                                                                   &\|(y_1,y_2)\|\ge k_{s_\mathrm{min}}\\
		\sqrt{k_{s_\mathrm{min}}^2-y_1^2-y_2^2} - k_{p_\mathrm{min}} & \text{otherwise}
	\end{cases}
\]
\item SS mode: K-space is covered by the region between two internally tangent circles $(-k_{s_\mathrm{max}}\vece_3,k_{s_\mathrm{max}})$ and $(-k_{s_\mathrm{min}}\vece_3,k_{s_\mathrm{min}})$. It consists of all points $\y\in\R^3$ such that $\|(y_1,y_2)\| \le k_{s_\mathrm{max}}$ and
\[
	\sqrt{k_{s_\mathrm{max}}^2-y_1^2-y_2^2} - k_{s_\mathrm{max}}
	\ge y_3
	\ge 
	\begin{cases} 
		-\|(y_1,y_2)\|                                                                   &\|(y_1,y_2)\|\ge k_{s_\mathrm{min}}\\
		\sqrt{k_{s_\mathrm{min}}^2-y_1^2-y_2^2} - k_{s_\mathrm{min}} & \text{otherwise}
	\end{cases}
\]
\end{itemize}
The different cases are depicted in Figure~\ref{fig:freqdiver} for transmission imaging (TI) and reflection imaging (RI). 
\begin{remark}
One remarks that the reflection imaging offers better spatial resolution for each modes by dint of the larger k-space coverage, however, it neglects lower spatial frequencies. 
\end{remark}

\section{Appendix}
\label{appendix} 
In this section, we will set some details of the proof of Theorem \ref{four_diff_theo} which basically has the same lines as the scalar case given in~\cite{KirQueRitSchSet21} with some changes. Also, mathematical details about distributions, (partial) Fourier transforms and convolutions can be found in~\cite{KirQueRitSchSet21}. Before that, one needs to calculate the partial Fourier transform for the Green's tensor $\mathcal{F}_{1,2}\vecG$ in the distributional sense. 

The Green's tensor $\vecG$ of an infinite, isotropic and elastic medium is given by 
\begin{equation}\label{gree_equ}
 \vecG(\x,\x'):=\frac{1}{\mu k_s^2}\left(k_s^2\bm{I}_2+\nabla\nabla\right)\mathcal{G}_s(\x,\x')-\frac{1}{\lambda+2\mu}\frac{\nabla\nabla}{k_p^2}\mathcal{G}_p(\x,\x')
 \end{equation}  
where $\mathcal{G}_s(\x,\x')=\frac{e^{ik_s\|\x-\x'\|}}{4\pi\|\x-\x'\|}$ and $\mathcal{G}_p(\x,\x')=\frac{e^{ik_p\|\x-\x'\|}}{4\pi\|\x-\x'\|}$ respectively are the full-space fundamental solutions for the Helmholtz equation with wavenumbers $k_s$ and $k_p$ satisfying the outgoing Kupradze radiation conditions, Equation~\ref{eq:kuprad}. 

Note that $\nabla\cdot\vecG_s=\bm{0}$ and $\nabla\times\vecG_p=\bm{0}$, where 
\[
\vecG_s(\x,\x'):=\frac{1}{\mu k_s^2}\left(k_s^2\bm{I}_2+\nabla\nabla\right)\mathcal{G}_s(\x,\x')\;\;\text{and}\;\;\vecG_p(\x,\x'):=-\frac{1}{\lambda+2\mu}\frac{\nabla\nabla}{k_p^2}\mathcal{G}_p(\x,\x')
\]
for more details see Equation 5.16.21 in~\cite{EriSuhCha78}. Afterwards, one defines 
\[
 \vecG_{\epsilon}(\x,\x'):=\frac{1}{\mu k_{s,\epsilon}^2}\left(k_{s,\epsilon}^2\bm{I}_2+\nabla\nabla\right)\mathcal{G}_{s,\epsilon}(\x,\x')-\frac{1}{\lambda+2\mu}\frac{\nabla\nabla}{k_{p,\epsilon}^2}\mathcal{G}_{p,\epsilon}(\x,\x')
\]
where $k_{{\alpha},\epsilon}=k_{\alpha}+i\epsilon$ the wavenumber and $\mathcal{G}_{\alpha,\epsilon}(\x,\x')=e^{-\epsilon\|\x-\x'\|}\mathcal{G}_{\alpha}(\x,\x')$, for $\epsilon>0$ the Green's function. Then, $\vecG_{\epsilon}\to\vecG$ in $\left[\it{S}'(\R^3)\right]^{3\times3}$ when $\epsilon\to 0$. 

One denotes the principal square root of $k_{\alpha,\epsilon}^2-\xi_1^2-\xi_2^2$ (i.e. the root with positive imaginary part) as follows 
\[
\kappa_{\alpha,\epsilon}:=\sqrt{k_{\alpha,\epsilon}^2-\xi_1^2-\xi_2^2},\quad\quad\text{for }\alpha=s,p
\]
\begin{lemma}\label{four_trans_gree}
The partial Fourier transform $\mathcal{F}_{1,2}\vecG\in\left[\it{S}'(\R^3)\right]^{3\times3}$ is given by 
\[
 \mathcal{F}_{1,2}\vecG\left[\varphibo\right]=\lim_{\epsilon\to 0}\mathcal{F}_{1,2}\vecG_{\epsilon}\left[\varphibo\right]
                                                                     =\lim_{\epsilon\to 0}\int_{\R^3}\left[ \hat{\vecG}_{s,\epsilon}e^{i\kappa_{s,\epsilon}|x_3|}- \hat{\vecG}_{p,\epsilon}e^{i\kappa_{p,\epsilon}|x_3|}\right]\varphibo(\xibo,x_3)\,\mbox{d}(\xibo,x_3)
\]
where
\[
 \hat{\vecG}_{s,\epsilon}=\frac{i}{4\pi\mu k_{s,\epsilon}^2\kappa_{s,\epsilon}}\left[k_{s,\epsilon}^2\bm{I}_2+\vecq_{s,\epsilon}\otimes\vecq_{s,\epsilon}\right],\quad\quad \hat{\vecG}_{p,\epsilon}=\frac{i}{4\pi\mu k_{s,\epsilon}^2\kappa_{p,\epsilon}}\vecq_{p,\epsilon}\otimes\vecq_{p,\epsilon}
\]
and 
\[
\vecq_{s,\epsilon}(\xibo)=-i\xibo'+i\kappa_{s,\epsilon}(\xibo)\text{sign}(x_3)\vece_3,\quad\quad\vecq_{p,\epsilon}(\xibo)=-i\xibo'+i\kappa_{p,\epsilon}(\xibo)\text{sign}(x_3)\vece_3
\]
for all Schwartz functions $\varphibo\in\left[\it{S}'(\R^3)\right]^3$.
\end{lemma}
\begin{proof}
The partial Fourier transform of the Green's tensor has the following form, using Equation~\ref{gree_equ} and formula $\frac{1}{\lambda+2\mu}\frac{1}{k_p^2}=\frac{1}{\mu k_s^2}$, as follows 
\begin{equation*}
\begin{split}
\mathcal{F}_{1,2}\vecG_{\epsilon}(\xibo,x_3,\x')
  &=\frac{1}{2\pi\mu k_{s,\epsilon}^2}\left[k_{s,\epsilon}^2\bm{I}_2+\vecq_{s,\epsilon}\otimes\vecq_{s,\epsilon}\right]\int_{\R^2}e^{-i(\xi_1x_1+\xi_2x_2)}\mathcal{G}_{s,\epsilon}(\x,\x')\,\mbox{d}x_1\mbox{d}x_2
\\&-\frac{1}{2\pi\mu k_{s,\epsilon}^2}\vecq_{p,\epsilon}\otimes\vecq_{p,\epsilon}\int_{\R^2}e^{-i(\xi_1x_1+\xi_2x_2)}\mathcal{G}_{p,\epsilon}(\x,\x')\,\mbox{d}x_1\mbox{d}x_2
\\&=\frac{1}{\mu k_{s,\epsilon}^2}\left[\left(k_{s,\epsilon}^2\bm{I}_2+\vecq_{s,\epsilon}\otimes\vecq_{s,\epsilon}\right)\mathcal{F}_{1,2}\mathcal{G}_{s,\epsilon}(\xibo,x_3,\x')-\left(\vecq_{p,\epsilon}\otimes\vecq_{p,\epsilon}\right)\mathcal{F}_{1,2}\mathcal{G}_{p,\epsilon}(\xibo,x_3,\x')\right]
\end{split}
\end{equation*}
Since $\mathcal{F}_{1,2}$ is a continuous operator on $\it{S}'(\R^3)$ (Proposition 7.17 in~\cite{KirQueRitSchSet21}), we have $\mathcal{F}_{1,2}\mathcal{G}_{\alpha,\epsilon}\to\mathcal{F}_{1,2}\mathcal{G}_{\alpha}$ in $\it{S}'(\R^3)$ and $\|\vecq_{\alpha,\epsilon}\otimes\vecq_{\alpha,\epsilon}\|_{\infty}\to\|\vecq_{\alpha}\otimes\vecq_{\alpha}\|_{\infty}$ when $\epsilon\to 0$, for $\alpha=s,p$. Consequently, we have $\mathcal{F}_{1,2}\vecG_{\epsilon}\to\mathcal{F}_{1,2}\vecG$ in $\left[\it{S}'(\R^3)\right]^{3\times3}$. The partial Fourier transform of the Green's function is given (Lemma 7.23 in~\cite{KirQueRitSchSet21}) as follows 
\[ 
 \mathcal{F}_{1,2}\mathcal{G}_{\alpha,\epsilon}=\mathcal{F}_{3}^{-1}\mathcal{F}\mathcal{G}_{\alpha,\epsilon}=\frac{i}{4\pi}\frac{e^{i\kappa_{\alpha,\epsilon}|x_3|}}{\kappa_{\alpha,\epsilon}}
\]
where $\mathcal{F}_{3}$ is the partial Fourier transform in the third direction. Finally, one gets
\begin{equation}\label{eq13}
\mathcal{F}_{1,2}\vecG_{\epsilon}(\xibo,x_3,\x')=\frac{i}{4\pi\mu k_{s,\epsilon}^2}\left[\left(k_{s,\epsilon}^2\bm{I}_2+\vecq_{s,\epsilon}\otimes\vecq_{s,\epsilon}\right)\frac{e^{i\kappa_{s,\epsilon}|x_3|}}{\kappa_{s,\epsilon}}-\left(\vecq_{p,\epsilon}\otimes\vecq_{p,\epsilon}\right)\frac{e^{i\kappa_{p,\epsilon}|x_3|}}{\kappa_{p,\epsilon}}\right]
\end{equation}
\end{proof}
\begin{proof}
Let $\vecg\in\left[L^p(\R^3)\right]^3$, for $p>1$, with $\text{supp}(\vecg)\subset\mathcal{B}_r$. Then, by the Sobolev embedding and the density of $\mathcal{D}(\R^3)$ in $L^p(\R^3)$, for $p\in[1,\infty)$, we can find a sequence of functions $\vecg_n\in\left[\mathcal{D}(\R^3)\right]^3$ with $\text{supp}(\vecg_n)\subset\mathcal{B}_r$, such that $\vecg_n\to\vecg$ in $\left[L^q(\R^3)\right]^3$, for $q\in[1,p]$ if $p\in(1,\infty)$ and for $q\in[1,\infty)$ otherwise, as $n\to\infty$. For each $\vecg_n$, one considers 
\[
   \left\{\begin{aligned}
    & \sigmabo_n(\x,\omega)=\vecC^0(\x):\varepsilonbo\left[\vecu_n\right](\x,\omega)                      \\
    & \nablabo\cdot\sigmabo_n(\x,\omega)+\omega^2\rho^0(\x)\vecu_n(\x,\omega) = -\vecg_n(\x,\omega)           
    \end{aligned} \right.
\]
with the outgoing Kupradze radiation conditions, Equation~\ref{eq:kuprad}. The unique solution $\vecu_n$ is given by the convolution $\vecu_n=\vecg_n\ast\vecG$ (Theorem 2, Section 5.16 in~\cite{EriSuhCha78}). Then, the partial Fourier transform is written as
\[
\mathcal{F}_{1,2}\vecu_n=2\pi\left(\mathcal{F}_{1,2}\vecg_n\right)\overset{3}{\ast}\left(\mathcal{F}_{1,2}\vecG\right)
\]  
where $\overset{3}{\ast}$ denotes the partial convolution with respect to the third coordinate, for more details see~\cite{KirQueRitSchSet21}. Now, for every $\varphibo\in\left[\it{S}(\R^3)\right]^3$, it follows by continuity of partial convolutions on $\it{S}'(\R^3)$ that 
\[
\mathcal{F}_{1,2}\vecu_n\left[\varphibo\right]=2\pi\lim_{\epsilon\to 0}\int_{\R^3}\mathcal{F}_{1,2}\vecG_{\epsilon}\left[M_3\mathcal{F}_{1,2}\vecg_n\overset{3}{\ast}\varphibo\right]\,\mbox{d}(\xibo,x_3)
\]
where $M_j(\varphi)(\x)=\varphi(x_1,...,x_{j-1},-x_j,x_{j+1},...,x_n)$. By Fubini's theorem and Equation~\ref{eq13}, one gets 
\begin{equation*}
\begin{split}
\mathcal{F}_{1,2}\vecu_n\left[\varphibo\right]
  &=\frac{i}{2\mu}\lim_{\epsilon\to 0}\frac{1}{k_{s,\epsilon}^2}\int_{\R^3}\{\left(k_{s,\epsilon}^2\bm{I}_2+\vecq_{s,\epsilon}\otimes\vecq_{s,\epsilon}\right)\frac{\varphibo}{\kappa_{s,\epsilon}}\int_{\R}e^{i\kappa_{s,\epsilon}|x_3-z|}\mathcal{F}_{1,2}\vecg_n(\xibo,z,\omega)\,\mbox{d}z
\\&-\left(\vecq_{p,\epsilon}\otimes\vecq_{p,\epsilon}\right)\frac{\varphibo}{\kappa_{p,\epsilon}}\int_{\R}e^{i\kappa_{p,\epsilon}|x_3-z|}\mathcal{F}_{1,2}\vecg_n(\xibo,z,\omega)\,\mbox{d}z\}\mbox{d}(\xibo,x_3)
\end{split}
\end{equation*}
Using the fact that $\mathcal{F}_{1,2}\vecg_n\in[L^1(\R^3)]^3$, one obtains by Lebesgue's dominated convergence theorem 
\begin{equation}\label{eq715}
\lim_{\epsilon\to 0}\frac{\varphibo(\xibo,x_3)}{\kappa_{\alpha,\epsilon}(\xibo)}\int_{\R}e^{i\kappa_{\alpha,\epsilon}|x_3-z|}\mathcal{F}_{1,2}\vecg_n(\xibo,z,\omega)\,\mbox{d}z=\frac{\varphibo(\xibo,x_3)}{\kappa_{\alpha}(\xibo)}\int_{\R}e^{i\kappa_{\alpha}|x_3-z|}\mathcal{F}_{1,2}\vecg_n(\xibo,z,\omega)\,\mbox{d}z
\end{equation}
for all $(\xibo,x_3)\in\R^3$ such that $\xi_1+\xi_2\neq k_{\alpha}^2,$ for $\alpha=s,p$. Using the fact that $\kappa_{\alpha}$ is locally integrable and $|\kappa_{\alpha}|\leqslant|\kappa_{\alpha,\epsilon}|$, it follows that 
\begin{equation}\label{eq716}
\left|\left(\vecq_{\alpha,\epsilon}\otimes\vecq_{\alpha,\epsilon}\right)\frac{\varphibo}{\kappa_{\alpha,\epsilon}}\int_{\R}e^{i\kappa_{\alpha,\epsilon}|x_3-z|}\mathcal{F}_{1,2}\vecg_n(\xibo,z,\omega)\,\mbox{d}z\right|\leqslant\|\vecq_{\alpha}\otimes\vecq_{\alpha}\|_{\infty}\frac{|\varphibo|}{2\pi|\kappa_{\alpha}|}\|\vecg_n\|_{L^1}\in[L^1(\R^3)]^3
\end{equation}
Taking into account Equation~\ref{eq715} and Equation~\ref{eq716} and applying Lebesgue's dominated convergence theorem, one gets
\begin{equation*}
\begin{split}
\mathcal{F}_{1,2}\vecu_n\left[\varphibo\right]
  &=\frac{i}{2\mu}\frac{1}{k_s^2}\int_{\R^3}\{\left(k_s^2\bm{I}_2+\vecq_s\otimes\vecq_s\right)\frac{\varphibo}{\kappa_s}\int_{\R}e^{i\kappa_s|x_3-z|}\mathcal{F}_{1,2}\vecg_n(\xibo,z,\omega)\,\mbox{d}z
\\&-\left(\vecq_p\otimes\vecq_p\right)\frac{\varphibo}{\kappa_p}\int_{\R}e^{i\kappa_p|x_3-z|}\mathcal{F}_{1,2}\vecg_n(\xibo,z,\omega)\,\mbox{d}z\}\,\mbox{d}(\xibo,x_3)
\end{split}
\end{equation*}
Therefore, we have ~\cite{KirQueRitSchSet21} 
\begin{equation}\label{eq718}
\begin{split}
\mathcal{F}_{1,2}\vecu_n\left[\varphibo\right]
  &=\sqrt{\frac{\pi}{2}}\frac{i}{\mu k_s^2}\int_{\R^3}\left(k_s^2\bm{I}_2+\vecq_s\otimes\vecq_s\right)\frac{\varphibo}{\kappa_s}\left(e^{i\kappa_s x_3}\mathcal{F}((1-H_{x_3})\vecg_n)(\xibo,\kappa_s,\omega)+e^{-i\kappa_s x_3}\mathcal{F}(H_{x_3}\vecg_n)(\xibo,-\kappa_s,\omega)\right)  
\\&-\left(\vecq_p\otimes\vecq_p\right)\frac{\varphibo}{\kappa_p}\left(e^{i\kappa_p x_3}\mathcal{F}((1-H_{x_3})\vecg_n)(\xibo,\kappa_p,\omega)+e^{-i\kappa_p x_3}\mathcal{F}(H_{x_3}\vecg_n)(\xibo,-\kappa_p,\omega)\right)\,\mbox{d}(\xibo,x_3)
\end{split}
\end{equation}
Considering $n\to\infty$ in Equation~\ref{eq718} for the right-hand side. Then, we obtain the point-wise limit with the fact that $\vecg_n\to\vecg$ in $[L^1(\R^3)]^3$ and $\mathcal{F}:L^1(\R^3)\to C_0(\R^3)$
\begin{equation}\label{eq719}
\begin{split}
  &\lim_{n\to\infty}\left(\vecq_{\alpha}\otimes\vecq_{\alpha}\right)\frac{\varphibo}{\kappa_{\alpha}}\left(e^{i\kappa_{\alpha} x_3}\mathcal{F}((1-H_{x_3})\vecg_n)(\xibo,\kappa_{\alpha},\omega)+e^{-i\kappa_{\alpha}x_3}\mathcal{F}(H_{x_3}\vecg_n)(\xibo,-\kappa_{\alpha},\omega)\right)
\\&=\left(\vecq_{\alpha}\otimes\vecq_{\alpha}\right)\frac{\varphibo}{\kappa_{\alpha}}\left(e^{i\kappa_{\alpha}x_3}\mathcal{F}((1-H_{x_3})\vecg)(\xibo,\kappa_{\alpha},\omega)+e^{-i\kappa_{\alpha}x_3}\mathcal{F}(H_{x_3}\vecg)(\xibo,-\kappa_{\alpha},\omega)\right)  
\end{split}
\end{equation}
for $\xi_1^2+\xi_2^2\neq\kappa_{\alpha}^2,\;\alpha=s,p$. As $\vecg_n\to\vecg$ in $[L^1(\R^3)]^3$, it exists a constant $C>0$ and $N\in\N$ such that 
\begin{equation*}
\begin{split}
& \left|\left(\vecq_{\alpha}\otimes\vecq_{\alpha}\right)\frac{\varphibo}{\kappa_{\alpha}}\left(e^{i\kappa_{\alpha}x_3}\mathcal{F}((1-H_{x_3})\vecg_n)(\xibo,\kappa_{\alpha},\omega)+e^{-i\kappa_{\alpha}x_3}\mathcal{F}(H_{x_3}\vecg_n)(\xibo,-\kappa_{\alpha},\omega)\right)\right|
\\&\leqslant\|\vecq_{\alpha}\otimes\vecq_{\alpha}\|_{\infty}\frac{|\varphibo|}{\kappa_{\alpha}}\|\vecg\|_{L^1}\in L^1(\R^3)
\end{split}
\end{equation*}
for every $n\geq N$ and almost every $(\xibo,x_3)\in\R^3$ such that $\xi_1^2+\xi_2^2\neq\kappa_{\alpha}^2$ for $\alpha=s,p$. Using the last two equations and the Lebesgue dominated convergence theorem, one gets  
\begin{equation}\label{eq721}
\begin{split}
  &\lim_{n\to\infty}\int_{\R^3}\left(\vecq_{\alpha}\otimes\vecq_{\alpha}\right)\frac{\varphibo}{\kappa_{\alpha}}\left(e^{i\kappa_{\alpha}x_3}\mathcal{F}((1-H_{x_3})\vecg_n)(\xibo,\kappa_{\alpha},\omega)+e^{-i\kappa_{\alpha} x_3}\mathcal{F}(H_{x_3}\vecg_n)(\xibo,-\kappa_{\alpha},\omega)\right)
\\&=\left(\vecq_{\alpha}\otimes\vecq_{\alpha}\right)\frac{\varphibo}{\kappa_{\alpha}}\left(e^{i\kappa_{\alpha}x_3}\mathcal{F}((1-H_{x_3})\vecg)(\xibo,\kappa_{\alpha},\omega)+e^{-i\kappa_{\alpha}x_3}\mathcal{F}(H_{x_3}\vecg)(\xibo,-\kappa_{\alpha},\omega)\right)  
\end{split}
\end{equation}

Finally, one considers the convergence of the left-hand side in Equation~\ref{eq718}. 
Limiting absorption principle (Theorem 1.1 in~\cite{AntFolPerRuiCru12}) gives that the unique solution $\vecu$ of Equation~\ref{eq32} for $\vecg\in[L^{q_1}(\R^3)]^3$, satisfying the outgoing Kupradze radiation conditions, fulfills 
\begin{equation}\label{eq722}
\|\vecu\|_{L^{q_2}}\leqslant C(k)\|\vecg\|_{L^{q_1}}
\end{equation}
if $\frac{1}{q_1}+\frac{1}{q_2}=1$ with $\frac12\leqslant\frac{1}{q_1}-\frac{1}{q_2}\leqslant\frac23$. Thus, for every $p>1$, one can find $1<q_1<p$ such as $q_1\in\left[\frac65,\frac43\right]$. For $\vecg\in[L^p(\R^3)]^3$ where $p>1$, one gets that $\vecg\in[L^{q_1}(\R^3)]^3$ where $1<q_1<p$ by the Sobolev embedding. Thus, Equation~\ref{eq722} implies that $\vecu_n\to\vecu$ in $[L^{q_2}(\R^3)]^3$, since $\vecg_n\to\vecg$ in $[L^{q_1}(\R^3)]^3$. In particular, $\vecu_n\to\vecu$ in $[\mathcal{S}'(\R^3)]^3$. Then, the continuity of $\mathcal{F}_{1,2}$ on $\mathcal{S}'(\R^3)$ gives 
 \[
\mathcal{F}_{1,2}\vecu\left[\varphibo\right]=\lim_{n\to\infty}\mathcal{F}_{1,2}\vecu_n\left[\varphibo\right] 
 \]
 for all $\varphibo\in[\mathcal{S}(\R^3)]^3$ and by Equation~\ref{eq721} finally 
 \begin{equation*}
 \begin{split}
 \mathcal{F}_{1,2}\vecu\left[\varphibo\right]&=\sqrt{\frac{\pi}{2}}\frac{i}{\mu k_s^2}\int_{\R^3}\left(k_s^2\bm{I}_2+\vecq_s\otimes\vecq_s\right)\frac{\varphibo}{\kappa_s}\left(e^{i\kappa_s x_3}\mathcal{F}((1-H_{x_3})\vecg)(\xibo,\kappa_s,\omega)+e^{-i\kappa_s x_3}\mathcal{F}(H_{x_3}\vecg)(\xibo,-\kappa_s,\omega)\right)  
\\&-\left(\vecq_p\otimes\vecq_p\right)\frac{\varphibo}{\kappa_p}\left(e^{i\kappa_p x_3}\mathcal{F}((1-H_{x_3})\vecg)(\xibo,\kappa_p,\omega)+e^{-i\kappa_p x_3}\mathcal{F}(H_{x_3}\vecg)(\xibo,-\kappa_p,\omega)\right)\,\mbox{d}(\xibo,x_3)
 \end{split}
 \end{equation*}
 for all $\varphibo\in[\mathcal{D}'(\R^3)]^3$. Then, the assertion follows by applying the du Bois-Reymond lemma, see Lemma 3.2 in~\cite{Gru09}.
 \end{proof}

\subsection*{Acknowledgments} 
This research was funded in whole, or in part, by the Austrian Science
Fund (FWF) P 34981. For the purpose of open access, the author has applied
a CC BY public copyright licence to any Author Accepted Manuscript version arising
from this submission. Moreover, OS is supported by the Austrian Science Fund (FWF),
with SFB F68 "Tomography Across the Scales", project F6807-N36 (Tomography with Uncertainties). The financial support by the Austrian Federal Ministry for Digital and Economic Affairs, the National Foundation for Research, Technology and Development and the Christian Doppler Research Association is gratefully acknowledged.
  
\section*{References}
\renewcommand{\i}{\ii}
\printbibliography[heading=none]

\end{document}